 \documentstyle[graphicx,amscd,amssymb,verbatim,12pt,righttag,color]{amsart}
 \newtheorem{theorem}{Theorem}[section]
 \newtheorem{Def}[theorem]{Definition}
 \newtheorem{Prop}[theorem]{Proposition}
 \newtheorem{Lem}[theorem]{Lemma}
 \newtheorem{Cor}[theorem]{Corollary}
 
 \newtheorem{Example}[theorem]{Example}

 \numberwithin{equation}{section}
 \renewcommand{\rm}{\normalshape}
 
\begin{document}

\title  {Classification of tile digit sets as product-forms}

\date{}
\author{Chun-Kit Lai}
\address{Department of Mathematics and Statistics, McMaster University,
Hamilton, Ontario, L8S 4K1, Canada}
\email{cklai@@math.mcmaster.ca}
\author{Ka-Sing Lau}
\address{Department of Mathematics , The Chinese University of
Hong Kong , Hong Kong}
\email{kslau@@math.cuhk.edu.hk}
\author {Hui Rao}
\address {Department of Mathematics, Central China Normal University, Wuhan, People's Republic of China}
\email{hrao@@mail.ccnu.edu.cn}
\thanks {}

\date{\today}
\thanks{ The research is supported in part by the HKRGC Grant, a Focus Investment Schemes in CUHK, and the CNSF (no. 11171100). }
\maketitle

\begin{abstract}
Let $A$ be an expanding matrix on ${\Bbb R}^s$ with integral entries. A fundamental question in the fractal tiling theory is to understand the structure of the digit set ${\mathcal D}\subset{\Bbb Z}^s$ so that  the integral self-affine set $T(A,\mathcal D)$ is a translational tile on ${\Bbb R}^s$. In our previous paper, we classified such tile digit sets ${\mathcal D}\subset{\Bbb Z}$ by expressing the mask polynomial $P_{\mathcal D}$ into product of cyclotomic polynomials. In this paper, we first show that a tile digit set in ${\Bbb Z}^s$ must be an integer tile (i.e. ${\mathcal D}\oplus{\mathcal L} = {\Bbb Z}^s$ for some discrete set ${\mathcal L}$). This allows us to combine the technique of Coven and Meyerowitz  on integer tiling on ${\Bbb R}^1$  together with our previous results  to characterize explicitly all tile digit sets ${\mathcal D}\subset {\Bbb Z}$ with $A = p^{\alpha}q$ ($p, q$ distinct primes)  as {\it modulo product-form} of some order, an advance of the previously known results for $A = p^\alpha$ and $pq$.

   \end{abstract}


\medskip

   \section{\bf Introduction}

\medskip

Let $A$ be an $s\times s$  expanding matrix (i.e. all eigenvalues have moduli $> 1$) with integral entries and  $|\det A |= b$ is a positive integer. Let  ${\mathcal{D}} \subset {\Bbb Z}^s$  and  call it a \textit{digit set}.  It follows that there exists
a unique compact set $T := T(A,{\mathcal{D}})\subset {\Bbb R}^s$ satisfying the set-valued relation $AT = T+{\mathcal{D}}$.  Alternatively, $T$ can be expressed as a set of radix expansions with base $A$ and digits in ${\mathcal{D}}$:
\begin{equation}\label{1.1}
T = \{\sum_{k=1}^{\infty}A^{-k}d_k : d_k\in
{\mathcal{D}}\}.
\end{equation}
It is well  known that  when  $\#{\mathcal D} = |\det A| = b$ and $T$ has non-empty interior, then $T$ is a translational tile in ${\Bbb R}^s$ \cite{[B]}. We call such $T$  a {\it self-affine tile} (or self-similar tile if $A$ is a scaling multiple of an orthonormal matrix) and ${\mathcal D}$  a {\it tile digit set} with respect to $A$. These tiles are referred as {\it fractal tiles} because their boundaries are usually fractals. There is a large literature on this class of tiles, and the reader can refer to them for the various developments ([LW1-4], \cite{[GH]}, \cite{[GM]}, \cite{[SW]}, \cite{[HLR]}, \cite{[KL]}, \cite{[LL]}, \cite{[GY]}). Note that for ${\mathcal D}$ to be a tile digit set, $\#{\mathcal D} = |\det A| = b $ is necessary. Hence in our consideration of tile digit sets, we will make this assumption without explicitly mentioning. For the digit sets in ${\Bbb R}^1$, we also assume, without loss of generality,  that ${\mathcal D} \subset {\Bbb Z}^+, \ 0 \in {\mathcal D}$ and g.c.d$({\mathcal D}) =1$.

\medskip

 Our interest is on the fundamental question of {\it characterizing the tile digit sets  ${\mathcal D}$ (i.e., $T(A, {\mathcal D})$ is a tile) for a given matrix $A$}. This turns out to be a very challenging problem even in ${\Bbb R}^1$.
 So far the only known cases are  $A=b$ with $b= p^{\alpha}$, a prime power \cite{[LW3]},  or $b= pq$, a product of two distinct primes \cite{[LR]}; there are extensions to the higher dimensional case when $|\det (A)| = p$ is a prime (\cite{[LW3]}, \cite{[HL]}). In this paper, we will advance our knowledge to the case when $b= p^\alpha q$, based on the theory that was developed in our previous paper \cite{[LLR]} and the new techniques here.

 \bigskip

 As is known, the most basic tile digit set ${\mathcal D}\subset {\Bbb Z}$ is ${\mathcal D} \equiv {\mathcal E} (\hbox {mod}\  b)$ where $b= \#{\mathcal D}$ and ${\mathcal E}=\{0, 1, \cdots , b-1\}$ (i.e., ${\mathcal D}$ is a complete residue set modulo $b$).  According to the earlier study in \cite{[LW3],[LR]}, it was suggested  that a tile digit set should possess certain product-form structure, i.e., ${\mathcal D}$ is obtained by decomposing ${\mathcal E}$ into direct sums according to the factors of  $b (=\#{\mathcal D})$, together with certain modulation by $b$. In \cite{[LLR]}, the authors investigated this idea of {\it modulo product-forms} in detail, and discovered the new classes of {\it $k$-th order modulo product-forms} of tile digit sets (see Definitions \ref {th3.1}, \ref{th3.2}). The main tool they used is the cyclotomic polynomials in elementary number theory, from which such product-forms are formulated algebraically in terms of product of cyclotomic factors of $P_{\mathcal D}$, the {\it mask polynomial} of ${\mathcal D}$. The very interesting part was the introduction of a tree structure on the cyclotomic polynomials ($\Phi$-tree), and the following theorem was proved (see also Section 3).

\medskip
\begin{theorem} {\rm \cite{[LLR]}} \
 ${\mathcal D} \subset {\Bbb Z}$ is a tile digit set if and only if  there is a blocking ${\mathcal N}$ in the $\Phi$-tree such that
$$
P_{\mathcal D} (x) = \Big (\prod_{\Phi_d\in{\mathcal N}}\Phi_d(x)\Big ) Q(x).
$$
\end{theorem}

\medskip

By a {\it blocking} in a tree, we mean a  finite subset ${\mathcal N}$  of vertices in the tree such that any paths from the root must pass through one and only one vertex in the set ${\mathcal N}$.  The above product is called a {\it kernel polynomial}. If ${\mathcal D}$ is a modulo product-form, then the kernel polynomial plays the role of product-form (see Section 3).

\medskip

In this paper, we continue our investigation along this line. First we will relate the tile digit sets with another well-known class of integer sets in number theory.  A finite subset ${\mathcal{A}}$ in ${\Bbb Z}^s$ is called an {\it integer tile} if there exists ${\mathcal{L}}\subset{\Bbb Z}^s$ such that ${\mathcal{A}}\oplus{\mathcal{L}} = {\Bbb Z}^s$ (The set $A\oplus B = \{a+b: a\in A, b\in B\}$ with $\oplus$ means that all elements are distinct). It is obvious that ${\mathcal{A}}$ is an integer tile if and only if ${\mathcal{A}}+[0,1]^s$ is a translational tile of ${\Bbb R}^s$.  The study of integer tiling dates back to the 40's, when Hajos and de Bruijn studied the factorization of abelian groups, and used it to solve some conjectures on Minkowski's geometry of
numbers.  One can refer to [Sz] for some details on the development of this topic. Our first main theorem, which holds in any dimension, is

\medskip

\begin{theorem}\label{th1.1}
Let $A$ be an integral expanding matrix. Suppose ${\mathcal  D} \subset {\Bbb Z}^s$ is a tile digit set (i.e.,  $T(A,{\mathcal{D}})$ is a self-affine tile), then ${\mathcal D}$ is an integer tile.
\end{theorem}

\medskip

The proof makes use of the self-replicating tiling sets of $T$ \cite{[LW2]} to construct the  tiling set $\mathcal L$  so that ${\mathcal D} \oplus {\mathcal L} = {\Bbb Z}^s$ (see Theorem \ref{th2.4}).  For the one-dimensional case, the theorem enables us to use the classical factorization techniques of cyclic groups to study the structure of tile digit sets. In particular, the most important technique we  employ is the decomposition of the integer tiles ${\mathcal A}$ when $\#{\mathcal A} = p^{\alpha}q^{\beta}$ developed by Coven and Meyerowitz \cite{[CM]}. We will combine this decomposition method and Theorem 1.1 to give the explicit expression of the kernel polynomial in the $\Phi$-tree (Theorem \ref{th5.5}), and to conclude the following (Theorem \ref{th5.6}).

\medskip

\begin{theorem}\label{th1.5}
Let ${\mathcal{D}} \subset {\Bbb Z}$ be a digit set with $\# {\mathcal{D}}= p^\alpha q$ and $p,q$ are primes. Then ${\mathcal D}$ is a tile digit set if and only if it is a  $k^{th}$-order modulo product-form for some $k$.
\end{theorem}

\medskip

The proof of Theorem \ref{th1.5} relies heavily on the algebraic operations on the cyclotomic polynomials (Proposition \ref{th2.6}), which  offers a lot more flexibility in handling the mask polynomial $P_{\mathcal D}$ than the set ${\mathcal D}$ itself. 

\medskip

We organize our papers as follows: In Section 2, we prove Theorem \ref{th1.1} and recall some needed facts on integer tiles. In Section 3, we summarize the notions of modulo product-forms and  cyclotomic tree we developed in  \cite{[LLR]}. In Section 4 and 5, we give a detailed study for  the case $b=p^{\alpha}q^{\beta}$ and prove Theorem \ref{th1.5} for the case $p^2q$, the case when $b=p^{\alpha}q$ is similar and is outlined. In Section 6, we give some further remarks on the tile digit sets and the integer tiles; the modulo product-forms and the related spectral problem of self-affine tiles are also discussed.

\bigskip

\section{\bf Tile digit sets on ${\Bbb R}^s$}

 We use the affine pair  $(A, {\mathcal D})$ in ${\Bbb R}^s$  and the attractor  $T := T(A, {\mathcal D})$ in ${\Bbb R}^s$ as in Section 1. Let  ${\mathcal{D}}_{A,k}:= \{\sum_{j=0}^{k-1}A^jd_{i_j}: d_{i_j}\in{\mathcal{D}}\} $ and ${\mathcal{D}}_{A,\infty}:=\bigcup_{k=1}^{\infty}{\mathcal{D}}_{A,k}$.
 The following are some well-known equivalent conditions for $T$ to be a tile  \cite{[LW2]}.

 \medskip

 \begin{theorem} \label{th2.1} That $ T := T(A, \mathcal D)$ is a self-affine tile (i.e., $T^o \not = \emptyset$) is equivalent to either one of the following conditions:

\medskip
\ \ (i) $\mu(T)>0$ where $\mu$ is the Lebesgue measure;

\medskip

\ (ii) $\overline{T^o} = T$ and $\mu(\partial T)=0$;

\medskip
(iii) $\#D_{A,k}=b^k$ for all $k \geq 0$. (Here $b = |\det (A)|$.)
\end{theorem}

\medskip

In this section our first goal is to show that if $T(A, {\mathcal D})$ is a self-affine tile (i.e.,  ${\mathcal D}$ is a tile digit set with respective to $A$), then ${\mathcal D}$ must be an integer tile in ${\Bbb Z}^s$. We need  a couple of lemmas to tie up some translational properties of the self-affine tiles . Let $\|\cdot\|$ denote the Euclidean norm.

\medskip

\begin{Lem}\label{lem2.2}
let $
E_k =\{ \sum_{j=1}^{\infty}A^{-kj}z:\  z\in{\mathcal{D}}_{A,k}\}$ and let $E = \bigcup_{j=1}^\infty E_k$,   then  $E$ is dense in $T$.
\end{Lem}

\medskip

\noindent {\bf Proof.} \ For any $y\in T$, write $y = \sum_{j=1}^{\infty}A^{-j}d_j$ and let  $y_k = \sum_{j=1}^{k}A^{-j}d_j$.  Let $z_k = A^k y_k$, then $z_k \in{\mathcal{D}}_{A,k}$, $\sum_{j=1}^{\infty}A^{-kj}z_k \in E_k$, and
$$
\|\sum_{j=1}^{\infty}A^{-kj}z_k -y_k\| = \|\sum_{j=1}^{\infty}A^{-kj}y_k\|  \leq \frac{\|A^{-k}\|}{1-\|A^{-k}\|}\ \|y_k\|\ .
$$
In view of the expanding property of $A$  and $\{\|y_k\|\}_{k=1}^{\infty}$ is bounded, the above expression tends to $0$ as $k$ tends to $\infty$.   This shows that $E$ is dense in $T$.
\qquad $\Box$

\bigskip

It is easy to see from (\ref{1.1}) that for any  $\ell\geq 1$,  $T(A,{\mathcal{D}})$ = $T(A^\ell,{\mathcal{D}}_{A,\ell})$. Moreover, if  we let $\widetilde{\mathcal D}: = {\mathcal D}_{A,\ell}-z^*$ be a translation of the digits, then by (\ref{1.1}),  $T(A^\ell,\widetilde {\mathcal{D}})$ satisfies
 \begin{equation}\label{2.1}
 T(A^\ell, \widetilde {\mathcal{D}})= T(A^\ell,{\mathcal{D}}_{A,\ell})-\sum_{j=1}^{\infty}A^{-j\ell}z^*.
 \end{equation}
We can impose some properties on the tile by suitably translating the digits.

\medskip

\begin{Lem}\label{lem2.3}
Suppose $T(A,{\mathcal{D}})$ is a self-affine tile, then there exists $z^{\ast}\in{\mathcal{D}}_{A,\ell}$ for some $\ell\geq1$ such that for  $\widetilde{{\mathcal{D}}}: = {\mathcal{D}}_{A,\ell}-z^{\ast}$, the tile $ {\widetilde T} :=  T(A^\ell,\widetilde{{\mathcal{D}}})$ satisfies
 $$
 0\in {\widetilde T}^o\quad \hbox {and}  \quad \partial \widetilde{T}\cap {\Bbb Z}^s = \emptyset .
 $$
 ($\partial \widetilde{T}$ denotes the boundary of $\widetilde{T}$.)
\end{Lem}

\medskip

\noindent {\bf Proof.} It is known that there exists  ${\mathcal{D}}': = {\mathcal{D}}_{A,k}-z, \ z \in{\mathcal{D}}_{A,k}$ such that $0\in T(A^k,{\mathcal{D}}')^o$ (\cite[Theorem 1.2]{[LW2]}).
Hence we can assume without loss of generality that $ 0 \in T^o \ (= T(A, {\mathcal D})^o)$.

\vspace {0.2cm}

Note that  $\partial T\cap{\Bbb Z}^s$ is a finite set. We let
\begin{equation} \label {2.3}
 \eta = \min\{\mbox{dist}(w,\partial T):  w \in {\Bbb Z}^s \setminus (\partial T\cap{\Bbb Z}^s) \}\ (>0).
\end{equation}
Let  $0< \epsilon < \eta$ such that $B(0, \epsilon ) \subset T$. It follows from $\mu (\partial T) =0$ (Theorem \ref {2.1}(ii)) and  the density of  $E$  in $T$ (Lemma \ref {lem2.2}) that there exists $ u \in B(0, \epsilon /2)$ with $u = \sum_{j=1}^\infty A^{-kj}z^*,  \ z^* \in {\mathcal D}_{A, k}$,  and $v + u \not \in \partial T$ for each $v\in \partial T\cap{\Bbb Z}^s$. This, together with \eqref {2.3}, yields
\begin{equation}\label{2.4}
(\partial T-u)\cap{\Bbb Z}^s =\emptyset .
\end{equation}
Let $\widetilde{{\mathcal{D}}} = {\mathcal{D}}_{A, k}-z^*$ and let  ${\widetilde T} =  T(A^k,\widetilde{{\mathcal{D}}})$. Then by (\ref{2.1}),
$$
{\widetilde T} = T - \sum_{j=1}^\infty A^{-kj}z^* = T - u.
$$
Now by our choice of $u$ and $\epsilon$, we have $
0\in B(0, {\epsilon}/{2})-u \subset \widetilde T$.
Hence $0\in {\widetilde T}^o$. Note also that $\partial {\widetilde T} = \partial T - u$,  we have, by \eqref {2.4}, $\partial {\widetilde T} \cap{\Bbb Z}^s = \emptyset$.  \qquad $\Box$

\bigskip

For a self-affine tile $T$, there exists a \textit{self-replicating} tiling set ${\mathcal{J}}\subset {\Bbb Z}^s$ \cite{[LW2]}, i.e., there exists $k\geq1$  such that
\begin{equation} \label {2.5}
A^k{\mathcal{J}}\oplus{\mathcal{D}}_{A,k} = {\mathcal{J}}.
\end{equation}
 Basically if $0\in T^o$, we can take ${\mathcal J} = {\mathcal D}_{A, \infty}$ in (\ref {2.5}). If $0\in \partial T$,  we translate $T$ so that $0$ is in the interior of the translated tile, then choose the $\mathcal J$ accordingly. The direct sum in  \eqref {2.5} is easy to check.  We now prove the main results in this section. Recall a finite set ${\mathcal A} \subset {\Bbb Z}^s$ is an {\it integer tile} if there exists ${\mathcal L}$ such that ${\mathcal A} \oplus {\mathcal L} = {\Bbb Z}^s$.

 \bigskip

\begin{theorem}\label{th2.4}
Let  $T(A,{\mathcal{D}})$ be a self-affine tile. Then ${\mathcal D}$ tiles ${\Bbb Z}^s$, i.e., ${\mathcal{D}}\oplus{\mathcal L} = {\Bbb Z}^s$ for some $\mathcal L$ in ${\Bbb Z}^s$.
\end{theorem}

\medskip

\noindent {\bf Proof.}
Let $\widetilde{{\mathcal{D}}} = {\mathcal{D}}_{A,\ell}-z^{\ast}$ be a digit set  so chosen that the conclusion of Lemma \ref{lem2.3} holds. As $0 \in T(A^\ell, \widetilde {\mathcal D})^o$,  it follows that ${\mathcal J}: = \widetilde{{\mathcal{D}}}_{A^l,\infty}$ is a self-replicating tiling set in ${\Bbb Z}^s$. Hence,
\begin{equation}\label{2.7}
A^\ell{\mathcal J}\oplus\widetilde{{\mathcal{D}}} = {\mathcal  J}.
\end{equation}

Let ${ \widetilde T} = T(A^\ell,\widetilde{{\mathcal{D}}})$ and let ${\mathcal B} = {\Bbb Z}^s\cap {\widetilde T}$. We  claim that ${\mathcal  J} \oplus {\mathcal B}  ={\Bbb Z}^s$. First, for any $w \in{\Bbb Z}^s$, we have $w \in {\widetilde T} + t $ for some $t \in{\mathcal J}$. This means that $w -t \in {\mathcal B}$. Hence, $ {\mathcal J} + {\mathcal B} = {\Bbb Z}^s $. To show that the representation is unique, we let  $w = t+ z = t'+ z'$ where $t, t'\in{\mathcal J}$ and $z, z'\in {\mathcal B}$. By the tiling assumption, we must have $w \in \partial {\widetilde T} +t$. This shows that $w-t \in {\mathcal B}\cap \partial {\widetilde T}$. But this is impossible since $\partial \widetilde{T}\cap{\Bbb Z}^s=\emptyset$ by our choice of $\widetilde{{\mathcal{D}}}$.

Now, by adding ${\mathcal B}$ to both side of (\ref{2.7}), we have
$$
A^\ell{\mathcal J}\oplus\widetilde{{\mathcal{D}}}\oplus {\mathcal B} = {\mathcal J}\oplus {\mathcal B} ={\Bbb Z}^s.
$$
Let ${\mathcal J}' = A^\ell{\mathcal{J}}\oplus {\mathcal B}$.  As
$\widetilde{{\mathcal{D}}} = {\mathcal{D}}_{A,\ell}-z^{\ast}$, we have from the above,
${\mathcal{D}}_{A,\ell}\oplus({\mathcal J}'-z^{\ast})={\Bbb Z}^s$. This implies
$$
{\mathcal{D}}\oplus(A{\mathcal{D}}\oplus \cdots \oplus A^{\ell-1}{\mathcal{D}}\oplus({\mathcal{J}}'-z^{\ast}))={\Bbb Z}^s.
$$
The theorem follows by setting ${\mathcal L} = A{\mathcal{D}}\oplus\cdots\oplus A^{\ell-1}{\mathcal{D}}\oplus( {\mathcal J}'-z^{\ast})$.
\qquad $\Box$

\bigskip

The converse of Theorem \ref {th2.4} is false in general, as is seen in the second part of the following example.

\medskip

\begin {Example} \label{ex2.5} \ Let $A=4$,  ${\mathcal D} = \{0,1,8,9\} = \{0,1\} \oplus 4\{0, 2\}$, then  $T(4, {\mathcal D}) =[0,1] \cup [2,3]$ is a self-similar tile (as it satisfies $4T= T+ {\mathcal D}$), and the tiling set for $T$  is  ${\mathcal J} = \{0, 1\} \oplus 4 {\Bbb Z}$. By Theorem \ref {th2.4},  ${\mathcal D}$ tiles ${\Bbb Z}$ also,  and the tiling set for ${\mathcal D}$ (as an integer tile) is  ${\mathcal L} = \{0, 2,  4, 6 \} \oplus 16{\Bbb Z}$.

\medskip

If we let ${\mathcal D} = \{0,1,4,5\} = \{0,1\} \oplus 4\{0, 1\}$, then ${\mathcal D}$ tiles ${\Bbb Z}$ with the tiling set  ${\mathcal L} = \{0,2\}\oplus 8{\Bbb Z}$. On the other hand, $T(4, {\mathcal D})$ is not a tile of ${\Bbb R}$ since $\#({\mathcal D} + 4{\mathcal D})  =12 < 4^2$  (by Theorem \ref{2.1}(iii)).
\end{Example}

\medskip

We remark that the explicit expression of the prime-power integer tiles and tile digit sets will be given in Theorems \ref{th4.2} and \ref {th5.1} respectively.

\bigskip

Throughout the rest of the paper, $a\mid b$ means $a$  divides $b$,   and $a\nmid b$ means $a$ does not divide $b$. The notations apply to both integers and polynomials. In the following we will give a brief summary on the cyclotomic polynomials and  integer tiles in ${\Bbb Z}^1$. Let $\Phi_d(x)$ be the \emph{$d$-th cyclotomic polynomial}, which is the minimal polynomial of the primitive $d$-th root of unity, i.e., $\Phi_d(e^{2\pi i/d}) = 0$. It is well-known that \
\begin{equation} \label {eq2.5'}
x^n -1 = \prod_{d|n}\Phi_d(x) \
\end{equation}
and the formula provides a constructive way to find $\Phi_d$ inductively.  The class of cyclotomic polynomials plays a fundamental role in the paper. Its basic manipulation rules are recalled below, they will be used extensively in Section 4 and 5.

\medskip

\begin{Prop}\label{th2.6}
Cyclotomic polynomials satisfy the following:

\vspace {0.2cm}

\ \  (i) If $p$ is a prime, then $\Phi_p(x) = 1+x+...+x^{p-1}$ and $\Phi_{p^{\alpha+1}}(x) = \Phi_p(x^{p^{\alpha}})$;

\vspace {0.2cm}

\ (ii) $\Phi_{s}(x^p) = \Phi_{sp}(x)$ if $p$ is prime and $p|s$, and

 \hspace {0.6cm} \ $\Phi_{s}(x^p) = \Phi_{s}(x)\Phi_{sp}(x)$ if $p$ is
prime but $p\nmid s$;

\vspace {0.2cm}

(iii) $\Phi_s(1) = \left\{
\begin{array}{ll}
0, & \hbox{if $s=1$}; \\
 p, & \hbox{if $s=p^{\alpha}$}; \\
  1, & \hbox{otherwise.}
 \end{array}
 \right.
 $
\end{Prop}

\medskip

    The class of integer tiles on ${\Bbb Z}$ has been studied in depth in connection with the factorization of cyclic groups and cyclotomic polynomials (\cite{[deB]}, \cite{[Tij]}, \cite{[CM]}, \cite{[S]},  \cite{[LW3]}, \cite{[N]}).  Let ${\Bbb Z}^+$ be the set of non-negative integers. For ${\mathcal{A}}\subset {\Bbb Z}^+$,  we let
  $$
  P_{\mathcal{A}}(x) = {\sum}_{a\in{\mathcal{A}}}x^a
  $$
   and call it the {\it mask polynomial} of ${\mathcal A}$.  The next simple lemma is well-known (see for example \cite{[CM]}) and it connects  cyclotomic polynomials with the factorization of ${\Bbb Z}_n$, the cyclic group with $n$ elements.

\medskip

\begin{Lem}\label{th2.7}
Let $n$ be a positive integer and let $\mathcal{A,B}$ be two finite sets of non-negative integers, then the following are equivalent.

\ \ (i) \ ${\mathcal{A}}\oplus{\mathcal{B}} \equiv {\Bbb Z}_{n}$,

\ (ii) \ $P_{\mathcal{A}}(x)P_{\mathcal{B}}(x) \equiv 1+x+...+x^{n-1}$ (mod $x^n-1$).

(iii) \ $n=P_{\mathcal{A}}(1)P_{\mathcal{B}}(1)$, and for every $d|n$, $\Phi_{d}(x)$ divides either $P_{\mathcal{A}}(x) $ or $P_{\mathcal{B}}(x)$.
\end{Lem}

\medskip

The notion in (i) means that ${\mathcal{A}}\oplus{\mathcal{B}} \equiv \{0, 1, \cdots, n-1\}$ (mod $n$). It is clear in this case   $\mathcal A$  is an {\it integer tile} with the {\it tiling set} ${\mathcal L} = {\mathcal B }\oplus n {\Bbb Z}$. (iii) relates the zeros $\{e^{2 \pi i/d}: \ d|n \} $ of $P_{\mathcal A}$ and $P_{\mathcal B}$ on the unit circle.

\medskip

For a finite set ${\mathcal A} \subset {\Bbb Z}^+$, we use
\begin{equation}\label{3.2}
S_{\mathcal{A}} = \{ p^\alpha > 1 : \ p \ \hbox {prime},
 \ \Phi_{p^\alpha}(x)|P_{\mathcal{A}}(x)\}
\end{equation}
to denote the {\it prime-power spectrum} of ${\mathcal A}$, and $\widetilde S_{\mathcal A} = \{s>1: \ \Phi_s(x)|P_{\mathcal{A}}(x)\}$ the {\it spectrum} of ${\mathcal A}$.  In \cite{[CM]}, Coven and Meyerowitz  made use of  the following two conditions to study the integer tiles:

\vspace {0.2cm}
\noindent {\bf (T1)} \  { \it $\#{\mathcal{A}} = P_{\mathcal{A}}(1)=\prod_{s\in
S_{\mathcal{A}}} \Phi_s(1)$,}

\vspace {0.1cm}

\noindent {\bf (T2)} \ {\it For any distinct prime powers $s_1, \ldots, s_n \in S_{\mathcal A} $ , then
$ s=s_1 \cdots s_n \in {\widetilde S}_{\mathcal{A}}$. }

\medskip

\begin{theorem} {\rm  \cite{[CM]}} \ \label{th2.8}  Let  ${\mathcal{A}}\subset {\Bbb Z}^+$ be a finite set, and satisfies  conditions $(T1)$ and $(T2)$, then ${\mathcal{A}}$ tiles ${\Bbb Z}$ with period \ $n =l.c.m.(S_{\mathcal{A}})$.

 Conversely, if $\mathcal{A}$ is an integer tile, then
$(T1)$ holds; if in addition  $\# {\mathcal A} = p^\alpha q^\beta, \ \alpha , \beta \geq 0$, then $(T2)$ holds.
\end{theorem}
It is still an open question whether an integer tile must satisfy (T2) in general.

   \bigskip
   \bigskip

   \section{\bf Modulo product-forms and cyclotomic trees}

\medskip

In this section, we will recall some  basic results about modulo product-forms and  cyclotomic trees developed in \cite{[LLR]}. They will be used for the explicit characterization of the tile digit sets ${\mathcal D}$ in Section 5.

\bigskip

Let $b\geq 2$, the product-form digit set of $b$  is defined as
\begin{equation} \label{eq3.1}
{\mathcal{D}} = {\mathcal{E}}_0 \oplus b^{l_1} {\mathcal{E}}_1 \oplus \ldots \oplus b^{l_k}
{\mathcal{E}}_k
\end{equation}
where
$  {\mathcal{E}} ={\mathcal{E}}_0 \oplus{\mathcal{E}}_1 \oplus \ldots \oplus{\mathcal{E}}_k \equiv {\Bbb{Z}}_b $, and $0\leq l_1 \leq
l_2\leq \ldots \leq l_k$ \cite{[LW3]};  if ${\mathcal{E}} = \{0,1,2,\ldots ,b-1\}$, then ${\mathcal{D}}$ is called a strict product-form  \cite{[O]}. The tile digit set in Example \ref{ex2.5} is such an example. However, such simple expression is far from covering all tile digit sets even when $b=4$.  In \cite{[LLR]}, we have introduced some more general classes of product-forms as tile digit sets.

\medskip

First we consider the product-form (\ref{eq3.1}) in terms of the mask polynomial and the product of cyclotomic polynomials. Observe that
$$
P_{\mathcal E}(x)= P_{{\mathcal E}_0}(x) P_{{\mathcal E}_1}(x)\cdots P_{{\mathcal E}_k}(x) = \prod_{d|b, d>1}\Phi_d(x)Q(x),
$$
and
 $$
P_{\mathcal D}(x) = P_{{\mathcal E}_0}(x) P_{{\mathcal E}_1}(x^{b^{l_1}})\cdots P_{{\mathcal E}_k}(x^{b^{l_k}})= \prod_{d|b, d>1}\Phi_d(x^{b_d})Q'(x)
$$
where $b_d$ is defined in the obvious way.  Based on the above product, we will generate more tile digit sets by taking modulo on each component. To this end, We need more notations. Let $S_i =\{d >1: d|b, \ \Phi_d(x)|P_{{\mathcal{E}}_i}(x)\}$ and let
\begin{equation*}
\Psi_i(x) = \prod_{d\in S_i}\Phi_{d}(x).
\end{equation*}
Then  $\Psi_i(x)|P_{{\mathcal{E}}_i}(x)$, hence $\Psi_i(x^{b^{l_i}})|P_{{\mathcal{E}}_i}(x^{b^{l_i}})$. Let
\begin{equation}\label{eq3.2}
K_1^{(i)}(x) = \Psi_0(x)\Psi_1(x^{b^{l_1}})...\Psi_i(x^{b^{l_i}}), \quad 0 \leq i\leq k \ .
\end{equation}
It is clear that for the product form, $K_1^{(k)}(x)= \prod_{d|b, d>1}\Phi_d(x^{b_d})$, which divides $P_{{\mathcal D}}(x)$.

\medskip

\begin{Def}\label{th3.1} Let ${\mathcal E} = {\mathcal{E}}_0 \oplus \ldots  \oplus  {\mathcal{E}}_k \equiv {\Bbb Z}_b$ and $0\leq l_1  \ldots \leq l_k $. For \
$
n_i   =  l.c.m.\ \{s: \ \Phi_{s}(x)\ |\ K_1^{(i)}(x)\},
$
 we define ${\mathcal{D}}^{(0)} \equiv {\mathcal{E}}_0 \ ({\mbox{mod}} \ n_0)$ and
\begin{equation} \label {eq3.3}
{\mathcal{D}}^{(i)} \equiv {\mathcal{D}}^{(i-1)} \oplus b^{l_i}{\mathcal{E}}_i \ ({\mbox{mod}} \ n_i), \qquad  1\leq i \leq k.
\end{equation}
We call ${\mathcal D} = {\mathcal D}^{(k)}$ the modulo product-form  with respect to  ${\mathcal E}$.
\end{Def}

\bigskip

It is clear that  product-form is a special case of modulo product-form. We also note that in the  mask polynomial, (\ref{eq3.3}) is equivalent to
$$
P_{{\mathcal{D}}^{(i)}}(x) =P_{{\mathcal{D}}^{(i-1)}}(x)P_{{\mathcal{E}}_{i}}(x^{b^{l_i}})+(x^{n_i}-1)
Q_{i+1}(x).
$$
By the choice of $n_i$ in (\ref{eq3.3}), we can prove that $K^{(i)}(x)|P_{{\mathcal D}}(x)$. This is used to show that the  modulo product-form is a tile digit set.  This extension, however, still does not cover all tile digit sets, as was shown by an example of ${\mathcal D}$ for $b=12$ in \cite[Example 3.6]{[LLR]}.  On the other hand, the example suggests the following higher order product-forms. In view of this, we call the above ${\mathcal D}$ a {\it $1^{st}$-order modulo product-form}.

\medskip

 \begin{Def}\label{th3.2}
${\mathcal{D}}$ is
called a $2^{\rm nd}$-order product-form if it is a  product-form of $\mathcal G$ with $0\leq  r_1 \cdots \leq r_{\ell}$, where ${\mathcal G} ={\mathcal G}_0 \oplus{\mathcal G}_1 \oplus \ldots \oplus{\mathcal G}_{\ell}$,  and  ${\mathcal G} $ itself is a $1^{st}$-order modulo product-form (as in Definition \ref{th3.1}, possibly in another decomposition different from the ${\mathcal G}_i$).

\vspace {0.15cm}

For the above ${\mathcal G}$ (as a $1^{st}$-modulo product-form), we let $K_{\mathcal G} = K_1^{(k)}$ be as in \eqref{eq3.2},
$$
S^{(2)}_i= \{s: \Phi_s(x)|P_{{\mathcal G}_i}(x), \ \Phi_s(x)|K_{{\mathcal G}}(x) \},
$$
$K_2^{(i)}(x) = \prod_{j=0}^{i}\prod_{s\in S^{(2)}_j}\Phi_{s}(x^{b^{r_j}})$ and $n_i = l.c.m.\{s: \Phi_s(x)|K_2^{(i)}(x)\}$ for $0\leq i\leq \ell.$
We use the same procedure as in Definition \ref{th3.1} to define the $2^{nd}$-order modulo product-form.

\vspace{0.15cm}

Inductively, we can define the  $k^{th}$-order product-form and modulo product-forms.

\end{Def}

\bigskip

Roughly speaking, we can produce new tile digit sets as follows: we start with the basic tile digit set  ${\mathcal E}$ (complete residue class), we construct the $1^{st}$-order modulo product-forms. We then rearrange those digits to form a product, and use them to construct the $2^{nd}$-order  modulo product-forms,  and likewise for the higher orders.  The interesting question is whether these higher order modulo product-forms will characterize all the tile digit sets. For this, we reformulate the question by expressing the mask polynomial into cyclotomic polynomials, and make use of the algebraic operations in Proposition \ref{th2.6} to study the question.

\medskip

Let $b\geq 2$. We define a {\it tree of cyclotomic polynomials (with respect to b)},  which we  call it a \textit{$\Phi$-tree}, as follows: the set of vertices of this tree at level $1$ are $\Phi_d$, where $d|b$ and $d>1$;   the offsprings of $\Phi_{d'}$ in each level are the cyclotomic  factors of $\Phi_{d'}(x^b)$, they are determined by Proposition \ref {th2.6}(ii). All $\Phi_d$ in the tree are different \cite[Proposition 2.2]{[LLR]}, hence it has a well-defined tree structure (see e.g. Figure 1 in Section 5). We call a finite subset of vertices  ${\mathcal N}$ a \emph{blocking} if every infinite path starting from the root intersects exactly one element of ${\mathcal N}$. The following is one of the main results proved in \cite{[LLR]}.

\medskip

\begin{theorem} \label{th3.3}
Let ${\mathcal{D}}\subset {\Bbb Z}^+$ with $\#{\mathcal D}=b$. Then  ${\mathcal D}$ is a tile digit set (with respect to $b$)
if and only if there is a blocking ${\mathcal N}$ in the $\Phi$-tree such that
$$
P_{\mathcal D} (x) = \Big (\prod_{\Phi_d\in{\mathcal N}}\Phi_d(x)\Big ) Q(x).
$$

\end{theorem}

\medskip

Let us denote the above product  by $K(x)$ and call it  a {\it kernel polynomial} of ${\mathcal D}$. The kernel polynomials play a central role in Section 5. It is seen that for the $1^{st}$-order modulo product-form, $K(x) = K_1^{(k)}(x)$ of ${\mathcal E}$ in (\ref{eq3.2}) is its kernel polynomial; the  $K_2^{(k)}(x)$  in Definition \ref {th3.3} is also a kernel polynomial. It is possible that by varying $Q(x)$, the same kernel can represent different tile digit sets ${\mathcal D}$. On the other hand, there are kernel polynomials that do not generate any tile digit set (see Remark 3 of Theorem 5.6 in \cite{[LLR]}). For the case $b=p^{\alpha}q$,  we will determine the admissible kernel polynomials and show that all the tile digit sets are $k^{th}$-order modulo product-forms for some $k$. To do so, we need to first specify the structure of the prime-power spectrum of a tile digit set, which is to be presented in the following theorem \cite[Theorem 2.4]{[LLR]}.

\medskip

\begin{theorem}\label{th3.4}
 Let  $b= p_1^{\alpha_1} \cdots p_k^{\alpha_k}$  be the product of prime powers and let $ {\mathcal D}$ be a tile digit set of $b$. Then the prime power spectrum of ${\mathcal D}$ is given by $$
 S_{\mathcal D} = {\bigcup}_{j=1}^k \{p_j^a : a \in E_j\}
 $$
 where
 $
 E_j = \{a :  \Phi_{p^a_j}(x) \mid  P_{\mathcal{D}}(x)\}$  and $E_j\equiv \{0, \cdots, \alpha_j -1\} (mod \ \alpha_j)$.
\end{theorem}

\medskip

The theorem can be proved directly from the Kenyon's criterion of self-similar tiles \cite[Theorem 2.4]{[LLR]}.  It can also be seen from Theorem \ref{th3.3} by considering the branches of the prime-power factors $d$ of $b$ in the first level.

\medskip

\section{\bf The $p^{\alpha}q^{\beta}$ integer tiles}

\medskip

For integer tiles ${\mathcal A}$  in ${\Bbb Z}$, there are special structural results when $\#{\mathcal A} = p^\alpha q^\beta$, a product of two prime powers \cite{[CM]}. Our main purpose in this section is to apply these results to prove, among the other results,  a special factorization lemma for ${\mathcal A}$ of cardinality  $p^\alpha q$ (Lemma \ref{th4.7}), which will be essential in Section 5.  For convenience, we assume that ${\mathcal A }  \subset {\Bbb Z}^+$, also we can assume that  g.c.d.$({\mathcal A}) =1$ whenever it is needed \cite[Lemma 1.4(1)]{[CM]}. First,  We start with a general lemma.

 \medskip

\begin{Lem}\label{th4.2}
 Let $P(x)\in{\Bbb Z}^+[x]$ and suppose that $P(x) = (x^n-1) g(x) +h(x)$  where $g(x)$ and $h(x)$ are respectively the quotient and remainder of $P(x)$ when it is divided $(x^n-1)$. Then $g(x), h(x)\in{\Bbb Z}^+[x]$. Moreover, if the coefficients of $h(x)$ are $0$ or $1$ only, then the same is for the coefficients of $P(x)$.
\end{Lem}

\medskip

\noindent{\bf Proof.} Write $P(x) = \sum_{i=1}^{N}a_ix^{m_i}$, where $a_i>0$, and let $m_i=\ell_i+nr_i$ with $0\leq\ell_i<n$. Then
$$
\begin{aligned}
P(x)=& \sum_{i=1}^{N}a_i(x^{m_i}-x^{\ell_i}+x^{\ell_i})\\
    = &\sum_{i=1}^{N}a_ix^{\ell_i}(x^{nr_i}-1)+ \sum_{i=1}^{N}a_ix^{\ell_i} \\
    =& (x^n-1)\sum_{i=1}^{N}a_ix^{\ell_i}(1+x^n+...+x^{n(r_i-1)})
    +\sum_{i=1}^{N}a_ix^{\ell_i}.\\
\end{aligned}
$$
Hence, $g(x) = \sum_{i=1}^{N}a_ix^{\ell_i}(1+x^n+...+x^{n(r-1)})$  and $h(x) = \sum_{i=1}^{N}a_ix^{\ell_i}$ ($\deg h<n$). This means that $g(x), h(x)$ have non-negative coefficients.  The second part is also clear. Indeed, the condition on $h(x)$ implies all the $\ell_i$ are distinct, and the above expression of $P(x)$ implies that  its coefficients equal to some $a_i$, and hence either $0$ or $1$.
\qquad{$\Box$}

\bigskip

 The major techniques we use are two well-known decomposition theorems (Theorems \ref{th4.1} and \ref{th4.4}). The first one is  due to de Bruijn \cite{[deB]}. Let ${\Bbb Z}^{+}[x]$ denote the set of polynomials with non-negative integer coefficients, and $F(x)({\mbox{mod}} \ x^n - 1)$ the remainder of $F(x)$ divided by $ (x^n - 1)$.

\medskip

\begin{theorem}\label{th4.1} {\rm (de Bruijn)} Let $n=p^\lambda q^\mu$ where $\lambda,\mu\geq 0$.
 Suppose $f(x)\in {\mathbb Z}^+[x]$ and $\Phi_n(x)|f(x)$, then there exist polynomials $P(x), Q(x)\in {\mathbb Z}^+[x]$ such that
\begin {equation} \label {eq4.1}
f(x)({\mbox{mod}} \ x^n - 1)=P(x)\Phi_{p^{\lambda}}(x^{q^{\mu}})+Q(x)\Phi_{q^{\mu}}
(x^{p^{\lambda}}).
\end{equation}
\end{theorem}

\medskip
 The theorem allows us to give an explicitly characterization of all integer tiles of prime power as  modulo product-form. (Note that the modulo product-form is originally defined for tile digit sets in \eqref {eq3.3}, however the same modulo procedure also works for integer tiles.)

 \medskip

Given an integer tile ${\mathcal A} \subset {\Bbb Z}^+$ with $\#{\mathcal A} = p^{\alpha}$,  denote $S_{\mathcal A} = \{p^{k_1},...,p^{k_{\alpha}}\}$ where $1\leq k_1<...<k_{\alpha}$  (if  g.c.d.$({\mathcal A}) =1$, then $k_1 =1$, also condition $(T 1)$ ensures $S_{\mathcal A}$ has exactly $\alpha$ elements).  Set $ {\mathcal E}_{k_j-1}  = p^{k_j-1}\{0,1, \ldots, p-1\}$. For ${\mathcal A}' ={\mathcal{E}}_{k_1-1}\oplus...\oplus{\mathcal{E}}_{k_\alpha-1}$,
\begin{equation*}\label{8.0+}
P_{{\mathcal A}'}(x) = P_{{\mathcal E}_{k_1-1}}(x)...P_{{\mathcal E}_{k_\alpha-1}}(x) = \prod_{i=1}^{\alpha}\Phi_{p^{k_i}}(x).
\end{equation*}
It is seen that ${\mathcal A'}$ is a product-form integer tile (by Theorem \ref{th2.8}).

\medskip

\begin{theorem}\label{th4.3} Let ${\mathcal A} \subset {\Bbb Z}^+$ with $\#{\mathcal A} = p^{\alpha}$. Then ${\mathcal A}$ is an integer tile if and only if  ${\mathcal A}$ is a modulo product-form of  ${\mathcal{E}}_{k_1-1}\oplus...\oplus{\mathcal{E}}_{k_{\alpha-1}
-1}$, in the sense that,
${\mathcal A } = A^{(\alpha)}$ with ${\mathcal A}^{(0)} =\{0\}$ and
\begin{equation} \label{eq4.3}
{\mathcal{A}}^{(i)} \equiv {\mathcal{A}}^{(i-1)} \oplus {\mathcal{E}}_{k_i-1}\ ({\mbox{mod}} \ p^{k_i}),
\end{equation}
for some $1\leq k_1<...<k_{\alpha}, \ 1\leq i \leq \alpha$.
\end{theorem}

\medskip

\noindent {\bf Proof.}  If ${\mathcal A}$ is a modulo product-form as in (\ref {eq4.3}), then $P_{\mathcal A}(x) = \prod_{i=1}^{\alpha}\Phi_{p^{k_i}}(x) Q(x)$. Then it satisfies (T1) and (T2),  and Theorem \ref{th2.8} implies it is an integer tile.

\vspace {0.2cm}

Conversely, assume ${\mathcal A}$ is an integer tile, then $S_{{\mathcal A}} = \{p^{k_1},...,p^{k_{\alpha}}\}$. We observe that with $n= p^{k^\alpha}$, (\ref {eq4.3}) is reduced to
\begin{equation} \label {eq4.4}
P_{\mathcal{A}}(x) (\mbox{mod} \ x^{p^{k_{\alpha}}}-1) = \Phi_{p^{k_{\alpha}}}(x) Q_{\alpha-1}(x)
\end{equation}
 and $ Q_{\alpha-1}(x)$ has nonnegative coefficients as the remainder has non-negative coefficients by Lemma \ref{th4.2}  and deg $(Q_{\alpha -1})< p^{k_{\alpha}} - \hbox {deg}(\Phi_{p^{k_{\alpha}}}) = p^{k_{\alpha}-1}$. Clearly, $ Q_{\alpha-1}(1) = p^{\alpha-1}$.
Observe that
 $\Phi_{p^{k_{i}}}(x)|Q_{\alpha-1}(x)$ for each $i$, we can repeat the same argument to obtain
$$
Q_{\alpha-1 }(x) (\mbox{mod} \ x^{p^{k_{\alpha-1}}}-1) = \Phi_{p^{k_{\alpha}-1}}(x)\ Q_{\alpha-2}(x),
$$
where $Q_{\alpha-2}(x)$ has nonnegative coefficients and $Q_{\alpha-2}(1) = p^{\alpha-2}$. Inductively we reach $Q_1(x)$ with the following identity,
\begin{equation} \label {eq4.5}
Q_1(x)(\mbox{mod} \ x^{p^{k_{1}}}-1) =  \Phi_{p^{k_{1}}}(x)Q_0(x),
\end{equation}
where $ Q_{0}(x)$ has nonnegative coefficients and $Q_0(1)=1$. Since $ Q_{0}(x)$ has nonnegative coefficients, we must have $Q_0\equiv1$.

\medskip

We now claim that ${\mathcal A}$ must be a modulo product-form. Note that $Q_1(x)$ is a polynomial with coefficients $0$ or $1$ (by \eqref{eq4.5} and $Q_0(x) \equiv 1$),  $Q_1(x)$ determines a  digit set ${\mathcal A}^{(1)}$. By \eqref{eq4.5} again and observe that $\Phi_{p}(x^{p^{k_1-1}})= \Phi_{p^{k_1}}(x)$,
$$
{\mathcal A}^{(1)}\equiv {\mathcal E}_{k_1-1} \ \ (\mbox{mod} \ p^{k_1}).
$$
Now, from  $Q_2(x)(\mbox{mod} \ x^{p^{k_2}}-1) = \Phi_{p^{k_2}}(x) Q_1(x)$, we have $\deg Q_1<p^{k_2}-(p-1)p^{k_2-1} = p^{k_2-1} $. This means that  $\Phi_{p^{k_2}}(x) Q_1(x)$ is a polynomial with coefficients $0$ or $1$, and by Lemma \ref {th4.2}, the same is for $Q_2(x)$. Let ${\mathcal A}^{(2)}$  be the digit set determined by $Q_2$, we have
$$
{\mathcal A}^{(2)}\equiv {\mathcal A}^{(1)}\oplus {\mathcal E}_2 \ \ (\mbox{mod} \ p^{k_2}).
$$
 Continuing this process, we  finally reach ${\mathcal A}$ as in \eqref {eq4.4}. Hence ${\mathcal A}$ is a modulo product-form in \eqref{eq4.3}.
\qquad{$\Box$}

\bigskip

Theorem \ref{th4.3} provides a  characterization of the structure of integer tiles of prime powers as some modulo product-forms. In particular, the techniques used in the proof of the theorem will appear again in Section 5.  Next, we will state another well-known decomposition theorem for the $p^{\alpha}q^{\beta}$ integer tiles, which is also derived from  de Bruijn's theorem (\cite{[S],[CM]}).

\bigskip

\begin{theorem}\label{th4.4}
  If ${\mathcal A}$ is an integer tile with $\#{\mathcal A} = p^{\alpha}q^{\beta}$, where $p,q$ are primes and $\alpha,\beta\geq 1$, and ${\mathcal A}\oplus {\mathcal B}\equiv{\Bbb Z}_n$, then  there is a prime factor of $\#{\mathcal A}$, say $p$,  such that
  $$
  \hbox {either}  \quad {\mathcal A}\subset p{\Bbb Z} \quad \hbox{or}\quad {\mathcal B} \subset p{\Bbb Z}.
  $$
  In the latter case (e.g., {\rm g.c.d.}$({\mathcal A}) =1$), we have
\begin{equation}\label{eq4.6}
{\mathcal{A}} = \bigcup_{j=0}^{p-1}\big (\{a_j\}\oplus p{\mathcal{A}}_j\big ),
\end{equation}
where $a_j = \min\{a\in{\mathcal{A}}: a\equiv j ({\mbox {mod}} \ p)\}$, and
${\mathcal{A}}_j = \{n \geq 0:  a_j +np \in{\mathcal{A}}\}$
are integer tiles.  In this case $\{a_j: \ 0\leq j\leq p-1\}$ forms a complete residue set (mod $p$) and  all $\#{\mathcal{A}}_j$ are equal.
\end{theorem}

\bigskip

We call \eqref{eq4.6} a {\it decomposition of ${\mathcal A}$ along $p$}. Putting the above in terms of the prime-power spectrum, we have

\medskip

\begin{Cor}\label{th4.5}
Suppose ${\mathcal{A}}$ is an integer tile with $\#{\mathcal{A}} = p^{\alpha}q^{\beta}$ and g.c.d.$({\mathcal A}) =1$, then we have
$$
p \in S_{\mathcal A}\  \mbox{ or } \ q \in S_{\mathcal A}  \quad  \hbox {(can be both)}.
$$
 In the case $p \in S_{\mathcal A}$, ${\mathcal A}$ has a decomposition along $p$ as in (\ref{eq4.6}), and
\begin {equation} \label {eq4.7}
 S_{p{\mathcal{A}}_j} = S_{{\mathcal{A}}}\backslash\{p\} \qquad \forall \ \  0\leq j\leq p-1.
\end{equation}
\end{Cor}

\medskip

\noindent{\bf Proof.}
There is a finite set $\mathcal B$ such that $
{\mathcal{A}}\oplus{\mathcal{B}} \equiv{\Bbb Z}_n$.  Note that ${\mathcal A}$ cannot be a subset of ${p{\Bbb Z}}$ since g.c.d.$({\mathcal A}) =1$.  Theorem \ref{th4.4} implies that there is a prime factor, say $p$, such that
${\mathcal B}\subset{p{\Bbb Z}}$. Writing $P_{{\mathcal{B}}}(x) = P(x^p)$ and note that
 $$
 P_{{\mathcal{B}}}(e^{2\pi i /p}) = P((e^{2\pi i /p})^p) = P(1) = P_{{\mathcal{B}}}(1) = \#{\mathcal{B}}\neq 0,
 $$
which means  $\Phi_{p}(x)$ does not divide $P_{{\mathcal{B}}}(x)$. But $\Phi_p(x) | P_{\mathcal A}(x)P_{\mathcal B}(x)$ (Lemma \ref{th2.7}(ii)). This implies $\Phi_{p}(x)|P_{{\mathcal{A}}}(x)$, i.e., $p \in S_{\mathcal A}$.

\medskip
That ${\mathcal A}$ has decomposition along $p$ follows from (\ref {eq4.6}).
Notice that ${\mathcal A}\oplus {\mathcal B}\equiv {\mathbb Z}_{n}$ and
${\mathcal A}_j \oplus {{\mathcal B}/p} \equiv {\mathbb Z}_{n/p}\ .$
These imply that
$$
S_{{\mathcal A}}\cup S_{\mathcal B}=S_{{\mathbb Z}_n} \quad {and} \quad S_{p{\mathcal A}_j}\cup S_{\mathcal B}=S_{{\mathbb Z}_n}\setminus \{p\},
$$
and the last statement of the corollary follows.
 \qquad{$\Box$}

\bigskip

 In the following we use the above theorems to consider the cyclotomic factors of $P_{\mathcal A}$ with  $\#{\mathcal{A}} = p^{\alpha}q^\beta$.

 \medskip

\begin{Lem}\label{th4.6}
Let ${\mathcal{A}}$ be an integer tile  with  $\#{\mathcal{A}} = p^{\alpha}q^{\beta}$.  If ${\mathcal{A}}$ admits a decomposition along $p$,  and  $\Phi_{p^{\lambda}q^{\mu}}(x)|P_{{\mathcal{A}}}(x)$ with $\lambda, \mu \geq 1$, then we have

\vspace{0.15cm}
\ (i) If $\lambda \geq 2$, then
$\Phi_{p^{\lambda-1}q^{\mu}}(x)|P_{{\mathcal{A}}_j}(x)$ for all $ j=0,...,p-1$;

\vspace{0.15cm}
(ii)   $\Phi_{p^{\lambda}}(x)|P_{{\mathcal{A}}}(x)$ or
 $\Phi_{q^{\mu}}(x)|P_{{\mathcal{A}}}(x)$.
\end{Lem}

\medskip

\noindent {\bf Proof}.
 (i)  By assumption, we can write
\begin {equation} \label {eq4.8}
P_{\mathcal{A}}(x) = \Phi_{p^{\lambda}q^{\mu}}(x)\cdot Q(x) =\Phi_{p^{\lambda-1}q^{\mu}}(x^p)\cdot Q(x).
\end{equation}
We use  $Q_j(x), \ 0\leq j \leq p-1$ to denote the polynomials of  the terms $x^{j+kp}$ in $Q(x)$,   then $Q(x)=\sum_{j=0}^{p-1} Q_j(x)$.   We define $\widetilde Q_j(x)$ by  $Q_j(x)=x^j {\widetilde Q}_j(x^p)$. Together with the decomposition in  (\ref{eq4.6}), we have
$$
x^{a_0}P_{{\mathcal{A}}_0}(x^p)+...+x^{a_{p-1}}P_{{\mathcal{A}}_{p-1}}(x^p) = {\mathcal P}_{\mathcal A}(x)  = \sum_{j=0}^{p-1}x^j\Phi_{p^{\lambda-1}q^{\mu}}(x^p){\widetilde Q}_j(x^p).
$$
By comparing the terms of the two polynomials, which has sorted according to the residue class modulo $p$, we have
$$
x^{a_j}P_{{\mathcal{A}}_j}(x^p) = x^j\Phi_{p^{\lambda-1}q^{\mu}}(x^p){\widetilde Q}_j(x^p) \qquad \forall \ \ 0\leq j\leq p-1.
$$
This implies that  $\Phi_{p^{\lambda-1}q^{\mu}}(x)|P_{{\mathcal{A}}_j}(x)$.

\medskip

 (ii) We will prove the statement by using induction on $ k =\alpha +\beta $. Let $n = p^{\lambda}q^{\mu}$.
If $k=1$, say $\#{\mathcal{A}} =p$, then  (\ref{eq4.1})
implies that  $p = pP(1)+qQ(1)$, which implies that  $Q(1) = 0$, so that $Q(x)\equiv 0$. Therefore $\Phi_{p^{\lambda}}(x^{q^{\mu}})|P_{{\mathcal{A}}}(x)$ (by (4.1) again), and hence $\Phi_{p^{\lambda}}(x)|P_{{\mathcal{A}}}(x)$ (by Proposition \ref {th2.6}\ (ii)).

For the induction step, by factoring out the $\mbox{g.c.d.}$ of ${\mathcal A}$, we can assume that g.c.d.$({\mathcal A}) =1$. Let us assume that $\Phi_{p}(x)|P_{{\mathcal{A}}}(x)$ as in Corollary \ref{th4.5}, hence ${\mathcal A}$ has a decomposition along $p$ as in (\ref{th4.6}). If $\lambda=1$, then we are done. If $\lambda\geq 2$, then part (i) implies that $\Phi_{p^{\lambda-1}q^\mu}(x)|P_{{\mathcal{A}}_i}(x)$  for all $i = 0,...,p-1$.  As all $S_{{\mathcal{A}}_i}$ are identical   (by Corollary  \ref{th4.5}),  we have, by the induction hypothesis, $\Phi_{p^{\lambda-1}}(x)$ or $\Phi_{q^{\mu}}(x)$ must divide all $P_{{\mathcal{A}}_i}(x)$. It follows that $\Phi_{p^{\lambda}}(x)$ or $\Phi_{q^{\mu}}(x)$  divides $P_{{\mathcal{A}}}(x)$
\qquad $\Box$

\bigskip

The following strengthens Lemma \ref {th4.6}(ii) for the case $\#{\mathcal A} = p^\alpha q$.

\begin{Lem}\label{th4.7}
Let ${\mathcal{A}}$ be an  integer tile such that \ $\#{\mathcal{A}} = p^{\alpha}q$.  If \ $\Phi_{p^{\lambda}q^{\mu}}(x)|P_{{\mathcal{A}}}(x)$, then
 \begin{equation}\label {eq4.9}
  \Phi_{p^{\lambda}}(x)\nmid P_{{\mathcal{A}}}(x) \ \ \Rightarrow \ \
 \Phi_{q^{\mu}}(x^{p^\lambda})|P_{{\mathcal{A}}}(x).
 \end{equation}
\end{Lem}

\medskip

\noindent {\bf Proof}.  In view of de Bruijn's identity (Theorem \ref {th4.1}),
\begin{equation} \label {eq4.11}
P_{{\mathcal{A}}}(x)( \mbox{mod} \ x^{p^{\lambda}q^{\mu}} -1) = P(x)\Phi_{p^{\lambda}}(x^{q^{\mu}})+
Q(x)\Phi_{q^{\mu}}(x^{p^{\lambda}}),
\end{equation}
 it suffices to show that if  $\Phi_{p^{\lambda}}(x)\nmid P_{\mathcal{A}}(x)$ (i.e., $p^\lambda \not \in  S_{\mathcal A}$), then $P(x)\equiv 0$.  By Corollary \ref {th4.5}, we have $p\in S_{\mathcal A}$ or $q \in S_{\mathcal A}$. We will divide our proof into two cases.

\medskip

  \noindent {\bf Case 1}:  $q \in S_{\mathcal A} $.  Suppose $P(x)$  in (\ref{eq4.11}) is not zero, let $x^t$ be a term with positive coefficient.  By checking the terms of $x^t \Phi_{p^{\lambda}}(x^{q^\mu })$ in the  product $P(x)\Phi_{p^{\lambda}}(x^{q^\mu })$, and noting that the terms in (\ref{eq4.11}) are positive, we conclude that there exists   ${\mathcal C} \subset {\mathcal A}$ and
\begin{equation*}
 {\mathcal C}  \equiv  t + q^\mu \{0,p^{\lambda-1},...,(p-1)p^{\lambda-1}\}\ \ (\hbox{mod} \ p^\lambda q^\mu).
\end{equation*}
Since $\Phi_q(x)$ divides $P_{{\mathcal{A}}}(x)$,  ${\mathcal{A}}$ admits a decomposition along $q$:
$$
{\mathcal A} = \bigcup_{j=0}^{q-1}\big (\{a_j\}\oplus q{\mathcal A}_j \big ).
$$
As all elements  of  ${\mathcal C}$ are in the same residue class (mod $q$), we  have
 ${\mathcal C}\subset \{a_j\}\oplus q{\mathcal A}_j$ for the $j$ such that ${a_j} \equiv t\ (\mbox{mod} \ q)$. Hence
\begin{equation} \label{eq4.12}
{\mathcal C}':=q^{\mu-1}\{0,p^{\lambda-1},...,(p-1)p^{\lambda-1}\}\  \subset \  {\mathcal A}_j (\hbox{mod } {p^\lambda}).
\end{equation}
 Note that $\#{\mathcal A} = p^\alpha q$ implies $\#{\mathcal A}_j=p^\alpha$. Let $S_{{\mathcal A}_j}= \{p^{k_1},...,p^{k_\alpha}\}$,  then by assumption, $p^\lambda \not \in S_{{\mathcal A}_j}$. Let $\ell$ be the largest integer such that $k_\ell<\lambda$.  By Theorem \ref{th4.3}, ${\mathcal A}_j$ is of the modulo product-form defined by ${\mathcal A}_j^{(i)}, 1\leq i \leq \alpha$ as in (\ref {eq4.3}). By taking modulo of $p^{\lambda}$, we see that
$$
{\mathcal A}_j({\mbox{mod}} \  p^{\lambda}) = {\mathcal A}_j^{(\ell)} \ ({\mbox{mod}} \  p^{k_\ell}).
$$
This leads to a contradiction as ${\mathcal C}'$ is a set of $p$ distinct elements in ${\mathcal A}_j({\mbox{mod}} \  p^{\lambda})$, but it is a singleton set $\{0\}$ in $A_j^{(\ell)}({\mbox{mod}} \  p^{k_\ell})$. Hence $P(x)\equiv 0$.

\medskip

\noindent {\bf Case 2}:  $p \in S_{\mathcal A}$.  \ Notice that in this case, we must have $\lambda>1$. The proof of this case is by induction on $\alpha$ and making use of Case 1 in the induction step.

When $\alpha =1$, from (\ref{eq4.11}),  we have $pq=pP(1)+qQ(1)$. Therefore either $P(1)=0,  Q(1)=p$ \ or \ $ Q(1)=0, P(1)=q$, which implies either $P(x)\equiv 0$ or $Q(x)\equiv 0$.
 But $Q(x)\not \equiv 0$ (for otherwise, set $x=e^{2\pi i/p^\lambda}$ in (\ref{eq4.11}), then the left-hand side is not $0$ but the right-hand side is $0$, a contradiction).   Hence we must have  $P(x) \equiv 0$,  so that $\Phi_{q^{\mu}}(x^{p^{\lambda}})|P_{{\mathcal{A}}}(x)$.

\vspace {0.2cm}

Suppose the statement holds for $\alpha -1$.
Since $p\in S_{{\mathcal{A}}}$,  we have by (\ref {eq4.6}),
$
P_{\mathcal{A}}(x) = x^{a_0}P_{{\mathcal{A}}_0}(x^p)+...+x^{a_{p-1}}
P_{{\mathcal{A}}_{p-1}}(x^p).
$
 This, together with (\ref {eq4.7}) and  $\Phi_{p^\lambda}(x) \nmid P_{\mathcal A}(x)$, implies that
$$
\Phi_{p^{\lambda -1}}(x) \nmid P_{{\mathcal A}_i}(x) \qquad  \forall \ \  0\leq i\leq p-1.
 $$
On the other hand Lemma \ref{th4.6} (i) implies that  $\Phi_{p^{\lambda -1}q^{\mu}}(x)|P_{{\mathcal A}_i}(x)$ for all $0\leq i\leq p-1$. Now if $p \in S_{{\mathcal A}_i}$,  we apply the induction hypothesis to conclude that $\Phi_{q^{\mu}}(x^{p^{\lambda-1}})|P_{{\mathcal A}_i}(x)$ for all $i$;  if $q \in  S_{{\mathcal A}_i}$, we can draw the same conclusion by Case 1.  This implies that in either cases, $\Phi_{q^{\mu}}(x^{p^{\lambda}})|P_{{\mathcal A}}(x)$.
\qquad{$\Box$}

\bigskip

\section{\bf Tile digit sets for $b=p^{\alpha}q$}

\medskip

In this section we will give an explicit characterization of the tile digit sets for $b=p^{\alpha}q$ and show that they are  modulo product-form of some order.  We assume that ${\mathcal D}\subset {\Bbb Z}^+, \ 0\in {\mathcal D}$ and g.c.d.$({\mathcal D}) =1$ as before. First we consider the simple case $b = p^\alpha$. Since a tile digit set ${\mathcal D}$ is an integer tile (Theorem \ref {th2.4}), they must have the form in (\ref{eq4.3}). Furthermore by Theorem \ref{th3.4}, the spectrum  $S_{{\mathcal D}}=\{p^{k_1},...,p^{k_\alpha}\}$ (with ${k_1} =1$) satisfies, in addition,  that  $\{k_i = i + \alpha t_i\}_{i=1}^{\alpha}$ is a complete residue set modulo $\alpha$.
In this case, ${\mathcal E}_{k_i-1} = p^{\alpha t_i}{\mathcal E}_{i-1}$ where ${\mathcal E}_i = p^i \{0,1, \ldots, p-1\}$. Hence
$$
{\mathcal D}' :={\mathcal{E}}_{k_1-1}\oplus...\oplus{\mathcal{E}}_{k_{\alpha}-1} = {\mathcal E}_0\oplus...\oplus p^{\alpha t_{\alpha}}{\mathcal E}_{\alpha-1}.
$$
This is a product-form since ${\mathcal E}_0\oplus...\oplus {\mathcal E}_{\alpha-1}=\{0,1,...,p^{\alpha}-1\}$. The following theorem is immediate.

\bigskip

\begin{theorem}\label{th5.1}
Suppose $b= p^\alpha$ for some $\alpha \geq 1$ and ${\mathcal{D}}$ is a tile digit set of  size $b$, then ${\mathcal D}$ is a modulo product-form of ${\mathcal E} = {\mathcal E}_0\oplus...\oplus {\mathcal E}_{\alpha-1}=\{0,1,...,p^{\alpha}-1\}$.
\end{theorem}

\medskip

We remark that in \cite{[LW3]}, Lagarias and Wang have a characterization of  tile digit sets of prime power. Their expression is more complicate, but it is the same as the above in essence. For $b=pq$, it was proved in \cite{[LR]} and \cite{[LLR]} that

\medskip

\begin{theorem} \label {th5.2}
Let $b = pq$ and ${\mathcal D}$ be a tile digit set of $b$. Then the prime power spectrum is $S_{{\mathcal D}} = \{p, q^n\}$ or $\{p^n, q\}$.  In the first case
$$
{\mathcal D} = \big (\{0, 1, \cdots p-1\} (\hbox{mod}\ p) \big ) \oplus b^{n-1}\{0, p, \ldots , p(q-1)\} \ (\hbox{mod} \ b^n),
$$
 and the kernel polynomial is
$$
 K(x) = \Phi_p(x) \Phi_{q^n}(x^{p^{n}}).
$$
 \end{theorem}

\bigskip

 In the rest of this section, we make use of the cyclotomic polynomial techniques developed in the previous sections to give a study of the tile digit sets for $b= p^{\alpha}q$. As the situation is rather involved in notations, we only work out the case $b=p^2q$ in detail (Theorems \ref{th5.5}, \ref{th5.6}), and the same idea applies to $b= p^{\alpha}q$ (Theorems \ref{th5.8}, \ref{th5.9}).

\medskip

  Theorem \ref{th3.3} will be needed when characterizing the tile digit sets by the cyclotomic $\Phi$-tree.  Figure 1  illustrates the $\Phi$-tree for $p^2q$ in which $\Phi_{p^m}$ has two offsprings, $\Phi_{q^n}$ has three, all the $\Phi_{p^mq^n}$ have only one.  It follows that the descendants of $\Phi_{p^mq^n}$  only form one path with no branch, hence the blocking on this path is  rather restricted, this becomes essential in classifying the admissible blockings  (see the proof of Theorem \ref{th5.5}(i) and (ii)).

\begin{figure}[h]\label{Fig1}
\centerline{\includegraphics[width=8cm,height=3cm]{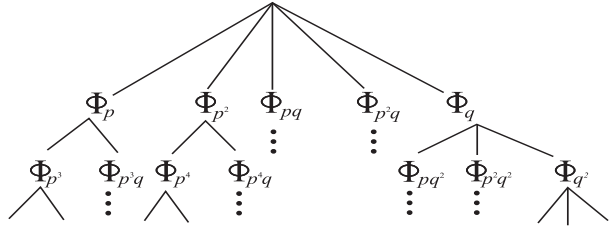}}
\caption {{\small $\Phi$-tree for  $p^2q$. }}
\end{figure}

\medskip

First we prove a simple lemma which will apply to the kernel polynomials.

\begin{Lem}\label{th5.4}
 Let ${\mathcal D}$ be a tile digit set (with g.c.d.$({\mathcal D})=1$ as assumed) and let
$G(x) = 1+ \sum_{j=1}^n a_jx^{k_j}, k_j \not =0$, be an integer polynomial
such that and $G(x) | P_{\mathcal D}(x)$ and $G(1)= \#{\mathcal D}$.  Then the
{g.c.d.} of $\{k_1,\dots, k_n\}$ is $1$.
\end{Lem}

\bigskip

\noindent{\bf Proof.} Suppose on the contrary that $\hbox {g.c.d.} \{k_1,\dots, k_n\}=d>1$. We can write $G(x) = \widetilde{G}(x^d)$ and $P_{\mathcal D}(x)= \widetilde{G}(x^d)Q(x) $. For each $0\leq j\leq d-1$, let $Q_j(x)$ be the polynomial containing the terms of the form  $x^{j+td}$ in $Q(x)$. Then $Q(x)=\sum_{j=0}^{d-1} Q_j(x)$, and we can also write $Q_j(x) =x^j\widetilde{Q}_{j}(x^d)$.  Hence
$$
P_{\mathcal D}(x) = \sum_{j=0}^{d-1} x^j \widetilde{G}(x^d)\widetilde{Q}_j(x^d).
$$
Let ${\mathcal D}_j = {\mathcal D}\cap (j+d{\Bbb Z})$, the subset of ${\mathcal D}$
that is congruent to $j (\hbox {mod}\ d)$. Then
$P_{\mathcal D}(x) = \sum_{j=0}^{d-1}P_{{\mathcal D}_j}(x)$.  By comparing the power, which has been sorted according to the residue
classes, we have for all $0\leq j\leq d-1$ that
 $$P_{{\mathcal D}_j}(x) = x^j \widetilde{G}(x^d)\widetilde{Q}_j(x^d).$$
  By taking $x=1$ and using the assumption that $G(1)( = \widetilde{G}(1)) = \#{\mathcal D}$, we see that $\#{\mathcal D}_j = \#{\mathcal D}\cdot Q_j(1)$. This means that $Q_j(1)$ is non-negative. Since $Q_j$ is a  monic polynomial of integer coefficients, $Q_j(1)$ is an integer.
As $\sum_j\#{\mathcal D}_j= \#{\mathcal D}$, it follows that only one of the $j$ satisfies $\#{\mathcal D}_j = \#{\mathcal D}$
and the others
are $0$. Since $0 \in {\mathcal D}_0$,
we must have ${\mathcal D}= {\mathcal D}_0$.  But then this contradicts to the fact that g.c.d.$({\mathcal D}) =1$ and hence the conclusion follows.
\qquad{$\Box$}

\bigskip

Let ${\mathcal{D}}$ be a tile digit set with $\#{\mathcal{D}} = p^2 q$ and g.c.d.$({\mathcal D}) =1$. Then it follows from Theorem \ref {th3.4} that
$$
S_{\mathcal{D}} = \{p, p^{2m}, q^n\}\  \mbox{ or }
 \ \{q, p^{2m}, p^{2n+1}\}.
$$
for some $m\in\{1,2,\cdots\}$ and $n\in \{0,1,\cdots,\}$. Our two main theorems are

\medskip

\begin{theorem}\label{th5.5}
Let $b = p^2q$ and let ${\mathcal{D}}$ be a tile digit set with $\#{\mathcal{D}} =b$. Then the mask polynomial $P_{\mathcal D}$ contains the following kernel polynomials:

\vspace {0.2cm}

\ {\bf (i)} \ If $S_{{\mathcal D}} = \{p,p^{2m}, q^n\}$,  then
\vspace {0.2cm}

\hspace {1.0cm} (I)\ \ \ \  $K_I(x) = \Phi_{p}(x)\Phi_{p^{2m}}(x^{q^{m}})\Phi_{q^n}(x^{p^{2(n-1)+1}})$; or
\vspace {0.2cm}

\hspace {0.9cm} (II) \ \  a factor of $ K_{II}(x) = \Phi_{p}(x)\Phi_{p^{2m}}(x^{q^{\ell-1}})
\Phi_{q^n}(x^{p^{2(n+m-\ell)}})$, $ 1\leq\ell\leq m$.

 \hspace {2.0cm}(as in (\ref{eq5.3}) below)

 \vspace {0.2cm}

 {\bf (ii)} \ If $S_{{\mathcal D}} = \{q,p^{2m}, p^{2n+1}\}$, then

 \vspace {0.2cm}

\hspace {0.8cm}  (III) \ \ $ K_{III}(x) = \Phi_{q}(x)\Phi_{p^{2m}}(x^{q^m})\Phi_{p^{2n+1}}(x^{q^{n+1}})$.

\vspace {0.2cm} Moreover, each of the above $K_i(x)$ represents a tile digit set of $p^2q$.

\end{theorem}

\bigskip

\begin{theorem}\label{th5.6}
Let $b = p^2q$ and assume that  $\#{\mathcal{D}} =b$. Then ${\mathcal D}$ is a tile digit set  if and only if it is a $k^{th}$-order modulo product-form for some $k\leq m$.

\end{theorem}

\bigskip

Using Proposition \ref{th2.6}(ii),  it is direct to check that  $K_{I}(x)$  is the mask polynomial of
\begin{equation} \label {eq5.1}
{\mathcal D}_I = {\mathcal E}_p \oplus b^{n-1}p{\mathcal E}_q \oplus b^{m-1}pq
{\mathcal E}_{p},
\end{equation}
where ${\mathcal E}_k = \{0, 1, \ldots , k-1\}$;  $K_{II}(x)$  is the mask polynomial of
\begin{equation}\label{eq5.2}
{\mathcal D}_{II} = {\mathcal E}_p \oplus b^{\ell-1}p^{2(m-\ell+1)-1}{\mathcal E}_p \oplus b^{n-1}p^{2(m-\ell+1)}{\mathcal E}_q,
\end{equation}
and $K_{III}$ is the mask polynomial of
\begin{equation*}
{\mathcal D}_{III} = {\mathcal E}_q \oplus b^n q{\mathcal E}_p \oplus b^{m-1} pq{\mathcal E}_{p}.
\end{equation*}

\medskip

\begin {Lem} \label{lem} With the above notations, ${\mathcal D}_I$ and ${\mathcal D}_{III}$ are $1^{st}$-order product-forms, and ${\mathcal D}_{II}$ is an $(m-\ell+1)$-order product-form.
\end{Lem}

\medskip

\noindent {\bf Proof.} Note that
$$
{\mathcal E}_p \oplus p{\mathcal E}_q \oplus pq{\mathcal E}_{p} = \{0, 1, \cdots, p^2q-1\} = {\mathcal E}.
$$
It follows that ${\mathcal D}_I$ is a product-form of ${\mathcal E}$. The proof for ${\mathcal D}_{III}$ is the same.

\vspace{0.2cm}

 For ${\mathcal D}_{II}$,  we let $t = m-\ell+1$, and let
 \begin{equation}\label {eq}
  {\mathcal D}' = {\mathcal E}_p \oplus p^{2t-1}{\mathcal E}_p \oplus p^{2t}{\mathcal E}_q.
 \end{equation}
 Then ${\mathcal D}_{II}$ is the product-form of ${\mathcal D}'$. It $t=1$, then clearly ${\mathcal D}_{II}$ is a $1^{st}$-order product form. Hence we  assume that $t>1$, observe that
 $
 {\mathcal E}_p \oplus p{\mathcal E}_q = q{\mathcal E}_p\oplus {\mathcal E}_q\  ( = \{0, 1, \cdots, pq-1\}),
 $
 we can rewrite ${\mathcal D}'$ as
 \begin{equation} \label {eq'}
 {\mathcal D}' = {\mathcal E}_p\oplus bp^{2t- 3}{\mathcal E}_p \oplus p^{2t-1}{\mathcal E}_q.
 \end{equation}
 Let ${\mathcal D}''= {\mathcal E}_p\oplus p^{2t-3}{\mathcal E}_p \oplus p^{2t-2}{\mathcal E}_q$.  We claim that ${\mathcal D}'$ is a modulo product form of ${\mathcal D}''$. Indeed, if we let $i\in{\mathcal E}_p$ and $j\in{\mathcal E}_q$, then $\{pj\}_{j\in{\mathcal E}_q}$ is a complete residue (mod $q$). This implies that
 $$
i+p^{2t-1}j  = i+p^{2t-2}(pj) = i+p^{2t-2}(qx_j+y_j) \equiv i+p^{2t-2}y_j \   (\mbox{mod} \ p^{2t-2}q)
 $$
and $\{y_j\}$ is a complete residue (mod $q$). Hence,
 $$
 {\mathcal E}_p\oplus p^{2t-1}{\mathcal E}_q\equiv {\mathcal E}_p\oplus p^{2t-2}{\mathcal E}_q \ (\mbox{mod} \ p^{2t-2}q)
 $$
and therefore ${\mathcal D}'$ is a modulo product-form of ${\mathcal D}''$ by ignoring the last modulo action.
   We continue this process for $t-1$ times  by noting that
    $$
   p^{2t-3}{\mathcal E}_p \oplus p^{2t-2}{\mathcal E}_q = p^{2t-3}({\mathcal E}_p\oplus p{\mathcal E}_q) =p^{2t-3}(q{\mathcal E}_p\oplus {\mathcal E}_q),
   $$
   and finally we obtain the first-order product form
$$
{\mathcal E} = {\mathcal E}_p\oplus p{\mathcal E}_p \oplus p^2{\mathcal E}_q \ (\equiv {\Bbb Z}_{p^2q}).
$$
This implies ${\mathcal D}_{II}$ must be an $(m-\ell+1)$-order product-form. \qquad  $\Box$

\bigskip

We remark that the rearrangement of (\ref{eq}) into (\ref{eq'}) is the key idea of the higher order modulo product-forms in Definition \ref{th3.2}. In term of cyclotomic polynomials, it  means that $\Phi_p(x)\Phi_{p^{2t}}(x)\Phi_q(x^{p^{2t}})= \Phi_p(x)\Phi_{p^{2t}}(x^q)\Phi_q(x^{p^{2t-1}})$ (by switching the position of $\Phi_q(x^{p^{2t}})$ from the last factor to the middle factor)  and the latter product is the modulo product-form of the  $\Phi_p(x)\Phi_{p^{2t-2}}(x)\Phi_q(x^{p^{2t-2}})$.

\bigskip

 For convenience, we call a vertex in a blocking ${\mathcal N}$ a {\it node}. Note that if  ${\mathcal D}$ is a tile digit set and $\Phi_d $ is a node in ${\mathcal N}_{\mathcal D}$, then $\Phi_d (x)| P_{\mathcal D}(x)$.  For the case $S_{\mathcal{D}} = \{p, p^{2m}, q^n\}$ in  Theorem \ref{th5.5}(i), the $K_{I}(x)$  gives a blocking of the $\Phi$-tree, and the nodes are determined by the following identities:
\begin{equation}\label{eq5.4}
 \Phi_{p^{2m}}(x^{q^{m}}) = \Phi_{p^{2m}}(x)\Phi_{p^{2m}q}(x)...\Phi_{p^{2m}q^{m}}(x)
\end{equation}
\begin{equation}\label{eq5.5}
\Phi_{q^n}(x^{p^{2(n-1)+1}}) = \Phi_{q^n}(x)\Phi_{pq^n}(x)...\Phi_{p^{2(n-1)+1}q^n}(x).
\end{equation}
Hence $K(x)$ in (I) is a kernel polynomial of ${\mathcal D}_I$.

\medskip

 $K_{II}(x)$  is a variant of $K_I$,  it is obtained by replacing the nodes $\Phi_{p^{2m}q^{i}}(x)$, $\ell\leq i\leq m$ (factors of $\Phi_{p^{2m}}(x^{q^m})$  in (I)) with new node $\Phi_{p^{2(m+n-i)}q^{n}}(x)$ (See Figure 1 and Lemma \ref{th5.7} below).
  Let
\begin{equation}\label{eq5.3}
{\widetilde K}_{II}(x)=\Phi_{p}(x)\Big ( \Phi_{p^{2m}}(x^{q^{\ell-1}})\prod_{j=n}^{n+m-\ell}
\Phi_{p^{2j}q^n}(x)\Big )\Phi_{q^n}(x^{p^{2(n-1)+1}}).
\end{equation}
It follows that ${\widetilde K}_{II}(x)$ is a kernel polynomial as its factors defines a blocking. Moreover $K_{II}(x)=  {\widetilde K}_{II}(x) {\widetilde Q}(x)$, with ${\widetilde Q}(x) = \prod_{j=n}^{n+m-\ell-1}
\Phi_{p^{2j+1}q^n}(x)$.

\bigskip

In order to classify the kernel polynomial of a tile digit set ${\mathcal D}$ in Theorem \ref{th5.5}, we need to know more precisely about the nodes. Let $\gamma_d$ denote the infinite path starts form $\vartheta$ and passes through $\Phi_d$. The following lemma describes the possible nodes on $\gamma_{p^\lambda q^\beta}, \lambda, \beta \geq 1$.

\medskip

\begin{Lem}\label{th5.7} With the above notations and $S_{\mathcal D} = \{p, p^{2m}, q^{n}\}$, then the nodes of ${\mathcal N}_{\mathcal D}$ satisfy

\vspace {0.2cm}

\ \ (i) for $0\leq k \leq m$, the node on $\gamma_{p^{2m}q^k}$ is either $\Phi_{p^{2m}q^k}$ or $\Phi_{p^{2(m+n-k)}q^n}$ ;

\vspace {0.2cm}

\ (ii) for $0\leq k \leq n-1$, the node on $\gamma_{{p^{2k+1}q^n}}$ is either $ \Phi_{p^{2k+1}q^n}$ or $\Phi_{pq^{n-k}}$ ;

\vspace {0.2cm}

(iii)  for  $1\leq k \leq n-1$, the node on $\gamma_{{p^{2k}q^n}}$ is either $\Phi_{p^{2k}q^n}$ or $\Phi_{p^{2m}q^{n+m-k}}$ .

\end{Lem}

\bigskip

\noindent{\bf Proof.}  To prove (i), note that the infinite path from $\Phi_{p^{2(m-k)}}$ pass through $\Phi_{p^{2m}q^k}$ has no other branch.  If $\Phi_{p^{2m}q^k}$ is not a node, then it must be an ancestor or a descendant of $\Phi_{{p^{2m}q^k}}$, i.e.,   $\Phi_{p^{2(m-r)}q^{k-r}}$ or $\Phi_{p^{2(m+r)}q^{k+r}}$. By Lemma \ref {th4.7} (or Lemma \ref{th4.6}(ii)), we must have $\Phi_{q^{k-r}}$ or $\Phi_{q^{k+r}}$ divides $P_{\mathcal D}$. They can only be $q^n$ since by assumption $S_{\mathcal A}= \{p, p^{2m}, q^n\}$. Hence $\Phi_{p^{2(m+n-k)}q^{n}}$ is the only choice. This completes the proof of (i).
The proof of (ii) and (iii) are similar. \qquad $\Box$

\medskip
Similarly, one can also develop an analogous lemma for $S_{\mathcal D} = \{q, p^{2m}, p^{2n+1}\}$. We can now prove our theorems.

\bigskip

\noindent {\bf Proof of Theorem \ref{th5.5}(i)}.  We divide the proof into two parts.

\vspace {0.2cm}

\noindent {\bf Case 1}:  Assume $\Phi_{p^{2m}}(x^{q^m})|P_{{\mathcal{D}}}(x)$, then $P_{\mathcal D}(x)= K_{I}(x)Q(x)$.

\medskip

 Note that if $n=1$, then $\Phi_{pq}(x) | P_{{\mathcal{D}}}(x)$ (by Lemma \ref {th5.7}(ii) with $k=0$). Together with $S_{{\mathcal D}} = \{p,p^{2m}, q\}$, we  have $\Phi_{p}(x)$, $\Phi_{q}(x)$, $\Phi_{pq}(x)$ and $\Phi_{p^{2m}}(x^{q^m})$ dividing $P_{{\mathcal{D}}}(x)$. Hence, $P_{\mathcal{D}}(x)$ must contain
$$
\Phi_{p}(x)\Phi_{p^{2m}}(x^{q^{m}})\Phi_{q}(x^{p}),
$$
which is of type (I).

\vspace {0.1cm}

For $n>1$, we claim that  $\Phi_{p^{2(n-1)+1}q^n}(x)| P_{{\mathcal{D}}}(x)$, then by observing that  $p^{2(n-1)+1}\not\in S_{{\mathcal{D}}}\ (= \{p, p^{2m}, q^n\}$) and applying Lemma \ref{th4.7},  $\Phi_{q^n}(x^{p^{2(n-1)+1}})$ will divide $P_{{\mathcal{D}}}(x)$,   and $P_{{\mathcal{D}}}(x)$ is of type (I).

Suppose otherwise,  $\Phi_{p^{2(n-1)+1}q^n}(x)\nmid P_{{\mathcal{D}}}(x)$. Then $\Phi_{pq}(x)| P_{{\mathcal{D}}}(x)$ (by Lemma \ref{th5.7}(ii)), hence $\Phi_p(x^q)| P_{{\mathcal{D}}}(x)$.
 Therefore $P_{\mathcal{D}}(x)$ must contain the following  factor
\begin{equation}\label{eq5.6}
G(x) = \Phi_{p}(x^q)\Phi_{p^{2m}}(x^{q^m})\Phi_{q^n}(x) = \Phi_{p}(x^q)\Phi_{p^{2m}}(x^{q^m})\Phi_{q^{n-1}}(x^q).
\end{equation}
It is now direct to check that $G(1) = p^2q = {\mathcal \#D}$, and  the g.c.d. of the non-zero power $k$ of $x^k$ in $G(x)$ is $q$.  This contradicts to Lemma \ref {th5.4} and hence Case 1 is proved.

 \medskip

\noindent {\bf  Case 2}:  Assume $\Phi_{p^{2m}}(x^{q^m}) \nmid P_{{\mathcal{D}}}(x)$, then $P_{\mathcal D}(x) = K_{II}(x)Q(x)$.

\vspace {0.2cm}

 Let $\ell$ be the first integer such that $1\leq \ell\leq m$ and $\Phi_{p^{2m}q^{\ell}}(x)\nmid P_{{\mathcal{D}}}(x)$, we claim that $\ell\neq n$. This is trivial if $m <n$. If  $n\leq m$, by taking $k=n$ in Lemma \ref{th5.7}(i), the two possibilities coincide as $\Phi_{p^{2m}q^n}(x)$. This means that $\Phi_{p^{2m}q^n}(x)$ must divide $P_{{\mathcal{D}}}(x)$ and hence $\ell \neq n$.

It follows that $m+n-\ell\neq m$.  By Lemma \ref{th5.7}(i), we have $\Phi_{p^{2(m+n-\ell)}q^n}(x)|P_{{\mathcal{D}}}(x)$. Since $p^{2(m+n-\ell)}\not\in S_{{\mathcal{D}}}$,  we must have $\Phi_{q^n}(x^{p^{2(m+n-\ell)}})|P_{{\mathcal{D}}}(x)$  (by Lemma \ref{th4.7}). Also the choice of $\ell$ implies that $\Phi_{p^{2m}}(x^{q^{\ell-1}})$ divides $P_{{\mathcal{D}}}(x)$. Hence $P_{{\mathcal{D}}}(x)$  contains a kernel polynomial of type (II).
\qquad$\Box$

\bigskip

\noindent {\bf Proof of Theorem \ref {th5.5}(ii)}.
In this case $S_{\mathcal{D}} = \{q, p^{2m}, p^{2n+1}\}$.  It follows from the (T2) property of integer tiles that
\begin {equation} \label {eq5.7}
\Phi_{p^{2m}q}(x)  |  P_{\mathcal D}(x) \ \ \hbox {and} \ \  \Phi_{p^{2n+1}q}(x)  |  P_{\mathcal D}(x)
\end{equation}
Also, analogous to Lemma \ref{th5.7}, we have

\vspace{0.2cm}
(i) \ for $1\leq k\leq m$, \ either \ $\Phi_{p^{2m}q^k}(x)|P_{{\mathcal{D}}}(x)$  \ or \ $\Phi_{p^{2(m-k+1)}q}(x)|P_{{\mathcal{D}}}(x)$; and

(ii) \ for  $1\leq l\leq n+1$, \ either \ $\Phi_{p^{2n+1}q^{l}}(x)|P_{{\mathcal{D}}}(x)$ \ or  \ $\Phi_{p^{2(n-l+1)+1}q}(x)|P_{{\mathcal{D}}}(x)$.

\vspace {0.1cm}
\noindent Let
$$
K_{III}(x) = \Phi_{q}(x)\Phi_{p^{2m}}(x^{q^m})\Phi_{p^{2n+1}}(x^{q^{n+1}}),
$$
then $K_{III}(x)$ is a kernel polynomial. We show that $P_{{\mathcal D}}(x)$ has $K_{III}(x)$ as a factor. Suppose otherwise, let $k_0$ and $l_0$ be the first integers such that
\begin {equation} \label {eq5.8}
\Phi_{p^{2m}q^{k_0}}(x)\nmid P_{{\mathcal{D}}}(x) \  \mbox{ and }  \  \Phi_{p^{2n+1}q^{l_0}}(x)\nmid P_{{\mathcal{D}}}(x) \quad \hbox {respectively}.
\end{equation}
  Note that $k_0>1$ (by (\ref {eq5.7})).  By (i), we have $\Phi_{p^{2(m-k_0+1)}q}(x)|P_{{\mathcal{D}}}(x)$.   Hence $p^{2(m-k_0+1)}\not \in S_{\mathcal{D}}$ (by (T2) of Theorem \ref{th2.8}).  It follows from Lemma \ref{th4.7} that $\Phi_{q}(x^{p^{2(m-k_0+1)}})|P_{{\mathcal{D}}}(x)$. By the same reasoning for the second part of \eqref {eq5.8},  we have $\Phi_q(x^{p^{2(n-l_0+1)+1}})|P_{{\mathcal{D}}}(x)$.  Let $\tau = \max\{2(m-k_0+1),2(n-l_0+1)+1\}>0$, we have $\Phi_{q}(x^{p^{\tau}})|P_{{\mathcal{D}}}(x)$. Hence, $P_{{\mathcal{D}}}(x)$ must contain the factor
$$
G(x) = \Phi_{q}(x^{p^{\tau}})\Phi_{p^{2m}}(x^{q^{k_0-1}})\Phi_{p^{2n+1}}
(x^{q^{l_0-1}})
$$
It is clear that $G(1) = p^2q$ and the $x^k$ of $G(x)$ has a common power $p$.  This contradicts  Lemma \ref{th5.4}, and hence, $K_{III}$ is a factor of $P_{\mathcal D}$.  This completes the proof.
\qquad{$\Box$}

\bigskip

\noindent{\bf Proof of Theorem \ref {th5.6}.}  The proof for the three types of $K(x)$ uses the same idea, we will only prove type (II) as it involves more variations. For this type,  there are two cases: (i) $n\geq\ell$ and (ii)$n<\ell$. For simplicity, we consider only the first case. The second case is similar by interchanging the last two factors in (\ref{eq5.9}) below. Let
\begin{equation}\label{eq5.9}
{\mathcal D}_{II} = {\mathcal{E}}_p \oplus b^{\ell-1}p^{2(m-\ell+1)-1}{\mathcal{E}}_p \oplus b^{n-1}p^{2(m-\ell+1)}{\mathcal{E}}_q.
\end{equation}
Then
\begin{equation} \label {eq5.8'}
K(x) := K_{II} (x) = \Phi_{p}(x)\Phi_{p^{2m}}(x^{q^{\ell-1}})
\Phi_{q^n}(x^{p^{2(n+m-\ell)}}).
\end{equation}
and the kernel polynomial is ${\widetilde K}(x)$  in (\ref{eq5.3}).
We will prove that ${\mathcal D}$ is a modulo product-form of ${\mathcal D}' = {\mathcal{E}}_p \oplus p^{2(m-\ell+1)-1}{\mathcal{E}}_p \oplus p^{2(m-\ell+1)}{\mathcal{E}}_q$. Then together with Lemma \ref{lem} (and also the proof),   ${\mathcal D}$ is an $(m-\ell+1)$-order modulo product-form.

\medskip

  To this end, we write $K(x) : = k_1(x)k_2(x)k_3(x)$ for the three factors in (\ref{eq5.8'}). We will use  the  similar technique as Theorem \ref{th4.3}.  Let
$n_3 \ = \ \hbox {l.c.m.}\{s: \Phi_{s}(x)|{\widetilde K}(x)\}$. As $\Phi_{p^{2(m+n-\ell)}q^n}(x)$ is in ${\widetilde K}(x)$. Then by definition $n_3 =  p^{2(m+n-\ell)}q^n$, and hence $k_3(x)|x^{n_3}-1$. Thus,
\begin{equation}\label{eq5.10}
P_{{\mathcal{D}}}(x)( \mbox{mod} \ x^{n_3}-1) =  K(x)Q(x).
\end{equation}
 Note that  ${ K}(x)Q(x)$ has non-negative coefficients (Lemma \ref{th4.2}), it implies that $k_1(x)k_2(x)Q(x)$ must also have non-negative coefficients. (In fact this follows from  $\deg (k_3) = \frac{n_3}{q}(q-1)$ and $k_3(x) = 1+x^{n_3/q}+ \cdots$, we have
$$
\deg(k_1k_2Q)<n_3-\deg k_3 = {n_3}/q,
$$
so that the terms of $k_1(x)k_2(x)Q(x)$ in the expansion of $ K(x)Q(x)$ do not overlap.)

\medskip

By considering $P_{\mathcal D}'(x) = k_1(x)k_2(x)Q(x)$ and letting
$n_{2} = \hbox {l.c.m.}\{s:\Phi_{s}(x)|k_1(x)k_2(x)\}$, we have $n_2 = p^{2m}q^{\ell}$ and
\begin{equation}\label{eq5.11}
{P}'_{\mathcal D}(x)(\mbox{mod} \ x^{n_2}-1) = k_2(x)(k_1(x)Q'(x)).
\end{equation}
By the same argument as the above,  $P_ {\mathcal D}''(x) = k_1(x)Q'(x)$ must have non-negative coefficients. Finally, let $n_1 =p$. We have
\begin{equation}\label{eq5.12}
P_{\mathcal D}''(x)(\mbox{mod} \ x^{p}-1)= \Phi_{p}(x)Q''(x).
\end{equation}
As $\deg (\Phi_{p}Q'')<p$ and $\deg(\Phi_{p})=p-1$, we must have $Q''(x)\equiv 1$. By combining (\ref{eq5.10}), (\ref{eq5.11}) and (\ref{eq5.12}), we see that ${\mathcal{D}} = {\mathcal D}^{(2)}$ where
$$
\left\{
  \begin{array}{ll}
    {\mathcal{D}}^{(0)} \equiv {\mathcal{E}}_p \ ({\mbox{mod}} \ n_1),\\
    \\
    {\mathcal{D}}^{(1)} \equiv {\mathcal{D}}^{(0)} \oplus b^{l-1}p^{2(m-l+1)-1}{\mathcal{E}}_p \ ({\mbox{mod}} \ n_2),\\
    \\
    {\mathcal{D}}^{(2)} \equiv {\mathcal{D}}^{(1)} \oplus b^{n-1}p^{2(m-l+1)}{\mathcal{E}}_q \ ({\mbox{mod}} \ n_3)
  \end{array}
\right.
$$
This proves the theorem.  \qquad{$\Box$}

\bigskip

For the case $p^\alpha q$, by Theorem \ref {th4.2} and  \ g.c.d.$({\mathcal D}) =1$, the prime power spectrum $S_{\mathcal D}$ is either

\vspace {0.2cm}

 {\bf (i)} $S_{\mathcal D} = \{p\}\cup \{p^{m_j \alpha +j}\}_{j=2}^{\alpha} \cup \{ q^n\}$ , or

\vspace {0.2cm}

{\bf (ii)} $S_{\mathcal D} = \{p^{m_j \alpha +j}\}_{j=1}^{\alpha} \cup \{q\}$

\vspace {0.2cm}

\noindent where $m_j \in {\Bbb N}^+$. (We modify the notation for $m_\alpha$ slightly in comparison with the $p^2q$ case, it is easy to check that $m_2+1=m$  for the  $m$ in Theorem \ref{th5.5}. The modification simplifies some expressions below.) We have

\medskip

\begin{theorem}\label{th5.8}
Let $b = p^{\alpha}q$ and let ${\mathcal{D}}$ be a tile digit set with $\#{\mathcal D} = b$. Then the mask polynomial $P_{\mathcal D}$ contains the following kernel polynomials

\vspace {0.2cm}

\noindent  Case {\bf (i)}:  either

\hspace {0.1cm}

\ (I) \ \ $  K_I(x) = \Phi_p(x) \Phi_{q^n}(x^{p^{\alpha(n-1)+1}}) \prod_{j=2}^{\alpha } \Phi_{p^{m_j\alpha+j}}(x^{q^{m_j+1}})$;  or
\vspace {0.1cm}

 (II) \ \ a factor of $ K_{II}(x) = \Phi_{p}(x)\Phi_{q^n}(x^{p^{{\alpha}(n+M)+k}})
\prod_{j=2 }^{\alpha } \Phi_{p^{m_j\alpha +j}}(x^{q^{\ell_j-1}})$,

\vspace {0.2cm}

\hspace {1.0cm}  where   \  $1\leq\ell_j \leq m_j+2  \ \forall \ j\geq 2$ and  at least one $\ell_j\leq m_j+1$ with

\hspace {1.0cm} $M = \max\{m_i-\ell_i: 2\leq i\leq \alpha \}$ and $ k = \max\{i: m_i-\ell_i = M\}$.

 \vspace {0.3cm}

\noindent Case {\bf (ii)}:

\vspace {0.2cm}
    (III) \ \ $ K_{III}(x) = \Phi_{q}(x) \prod_{j=1}^{\alpha }\Phi_{p^{ m_j\alpha +j}}(x^{q^{m_j+1}})$.

\vspace {0.3cm}
 Moreover, each of the above $K_i(x)$ represents a tile digit set of $b= p^\alpha q$.
\end{theorem}


\noindent{\bf Remark.} For the polynomial in $K_{II}(x)$, if $\ell_j =m_j+2$, then $\ell_j-1=m_j+1$ and this means that the whole factor $\Phi_{p^{m_j\alpha+j}}(x^{q^{m_j+1}})$ divides $P_{\mathcal D}(x)$. However, for a non-trivial $K_{II}(x)$ to occur, we must need at least one $\ell_j\leq m_j+1$. This is consistent with Theorem \ref{th5.5}, case (II), since there is only one factor in the above product when $b= p^2q$, which reduces to $1\leq \ell_2\leq m_2+1$.

\bigskip

It is direct to check that that for $K_{I}(x)$, it represents a tile digit set
$$
{\mathcal D}_I = {\mathcal E}_p\oplus b^{n-1}(p{\mathcal E}_q)\oplus \bigoplus_{j=2}^{\alpha} b^{m_j}(p^{j-1}q{\mathcal E}_p)
$$
 For $K_{II}(x)$,
$$
{\mathcal D}_{II} = {\mathcal E}_p\oplus b^{n-1}(p^{\alpha(M+1)+k}{\mathcal E}_q)\oplus \bigoplus_{j=2}^{\alpha} b^{\ell_j-1}(p^{(m_j-\ell_j+1)\alpha+j-1}{\mathcal E}_p).
$$
For $K_{III}(x)$,
$$
{\mathcal D}_{III} = {\mathcal E}_q\oplus \bigoplus_{j=1}^{\alpha} b^{m_j}(p^{j-1}q{\mathcal E}_p).
$$
These three digit sets are up to an rearrangement of the factors so that  the powers of $b$  are in non-decreasing order.

\bigskip

The proof is basically identical to the case  $b=p^2q$.  First, it is easy to deduce Lemma \ref{th5.7} for $b=p^{\alpha }q$.   For Case {\bf (i)} We first suppose that $ \prod_{j=2}^{\alpha} \Phi_{p^{m_j\alpha+j}}(x^{q^{m_j+1}})|P_{{\mathcal D}}(x)$. The same argument shows that $\Phi_{q^{n}p^{\alpha(n-1)+1}}(x)$ will divide $P_{\mathcal D}(x)$. Hence, we have type (I) using Lemma \ref{th4.7}. Suppose that the above factors does not divide $P_{\mathcal D}(x)$. For those $\Phi_{p^{m_j\alpha+j}}(x^{q^{m_j+1}})\nmid P_{{\mathcal D}}(x)$,   we let  $\ell_j$ be the first integer such that $\Phi_{p^{m_j\alpha+j}q^{\ell_j}}(x)$ does not divide $P_{{\mathcal D}}(x)$ (hence $\ell_j\leq m_j+1$).
By making use of the technique in  Theorem \ref{th5.5}(i)-Case(2), we obtain $\Phi_{q^{n}}(x^{p^{\alpha(n+m_j-\ell_j)+j}})|P_{{\mathcal D}}(x)$. Note that
$$
\max_{2\leq j\leq \alpha}\{\alpha(n+m_j-\ell_j)+j\} = \alpha(n+M)+k,
 $$
 where $M$ and $k$ are defined in the statement. Hence,  $\Phi_{q^{n}}(x^{p^{\alpha(n+M)+k}})$ divides $P_{\mathcal D}$(x). This forms $K(x)$ of type (II) which contains the kernel polynomial $K'(x)$. The previous consideration for type (III) applies in the same way.

\bigskip

\begin{theorem} \label{th5.9}
Let $b = p^{\alpha}q$ and let ${\mathcal{D}}$ be a tile digit set $\#{\mathcal D} = b$. Then ${\mathcal{D}}$ must be a $k^{th}$-order modulo product-form.
\end{theorem}

\medskip

It follows from the same idea of proof as Theorem \ref{th5.6}, The first thing is to establish the analogue of Lemma \ref{lem} for  ${\mathcal D}_I$, ${\mathcal D}_{II}$ and ${\mathcal D}_{III}$ given above. It is easy to see that ${\mathcal D}_I$ and ${\mathcal D}_{III}$ are 1-st order product-form. For ${\mathcal D}_{II}$, pick $j$ so that $m_j-\ell_j =M $ and $k=j$, then we note that
$$
p^{\alpha(M+1)+k}{\mathcal E}_q\oplus p^{(m_j-\ell_j+1)\alpha+j-1}{\mathcal E}_p = p^{(m_j-\ell_j+1)\alpha+j-1}({\mathcal E}_p\oplus p{\mathcal E}_q)=p^{(m_j-\ell_j+1)\alpha+j-1}(q{\mathcal E}_p\oplus {\mathcal E}_q).
$$
From this, we can use the same argument in Lemma \ref{lem} to conclude that  ${\mathcal D}_{II}$ is a product-form of some orders $t$.

\medskip

Next, we need to show that any tile digit sets must be given by the modulo product-forms of ${\mathcal D}_I$, ${\mathcal D}_{II}$ and ${\mathcal D}_{III}$.  This is done by arranging the powers of $b$ of those digit sets in non-decreasing order and applying the same argument in the proof of Theorem \ref{th5.6}. We can eventually show that all are of modulo product-form of order $t+1$.

\bigskip

\section{\bf Some remarks}

\medskip
 One of the aims in our investigation is to study the role of the modulo product-forms (and the higher order ones) in the tile digit sets, in particular, to characterize the tile digit sets to be such forms. So far we can only describe the modulo product-forms on ${\Bbb R}^1$. It will be interesting to define the analog in the higher dimensional spaces. Note that the definition of the product-form is easy to be generalized (see. e.g. \cite{[LW3]}), however, there is no direct generalization for the modulo product-form. One of the main difficulties is to find some replacement of the cyclotomic polynomials.
In another direction, it will also be useful to develop an algorithm to check for a given digit set ${\mathcal D}\subset {\Bbb Z}^+$ to be a tile digit set.

 \bigskip

  The main techniques we use in the explicit characterization for the tile digit sets of $\#{\mathcal D} = p^\alpha q$ is  the classical results  of de Bruijn about $\Phi_{ p^\alpha q^\beta}(x) | P_{\mathcal D}(x)$ (Theorem \ref{th4.2}), and the decomposition of integer tiles ${\mathcal A}$ when $\#{\mathcal A}= p^\alpha q^\beta$ (Theorem \ref{th4.4}). It is likely that our approach can  further be improved to obtain a complete characterization of tile tile digit sets of $\#{\mathcal D} = p^\alpha q^\beta$ (as well as integer tiles ${\mathcal A}$ of the same cardinality) as certain kind of modulo product-forms as in Theorems \ref{th5.8}, \ref{th5.9}. Finally and more challengingly, if $\#{\mathcal A}$ or $\#{\mathcal D}$ has more than two prime factors, some new factorization theorems may need to be developed.

\medskip

Our study of the  structure  of the tile digit sets is closely related to the spectral set problems. Recall that a closed subset $\Omega \subset {\Bbb R}^s$ is called a {\it spectral set} if $L^2(\Omega)$ admits an exponential orthonormal basis $\{e^{2\pi i \langle \lambda, \cdot \rangle}\}_{\lambda\in\Lambda}$ ($\Lambda$ is called a spectrum ). The well-known Fuglede conjecture asserted that $\Omega$ is a spectral set if and only if it is a translational tile. The conjecture was eventually proved to be false on ${\Bbb R}^s$ for $s\geq 3$ (\cite{[T],[KM]}).  The conjecture is still widely open for self-affine tiles. Our consideration on cyclotomic polynomials factors  for the tile digit sets is  closely linked to the spectral problem, because it also deals with zeros on the unit circle. It would be instructive to first study the spectral problem for simpler product-form as a testing case. On the other hand, there are studies of the spectral problem for integer tiles on ${\Bbb R}^1$ \cite{[L]}. In view of Theorem \ref{th2.4}, the results can be applied to the tile digit sets, and it may offer some insight to investigate the spectral  problem of the self-affine tiles.

   \end{document}